\documentclass{article}

\usepackage{fullpage, graphicx, amsmath, amssymb}
\usepackage[tight]{subfigure}
\usepackage[ruled, vlined, algo2e]{algorithm2e}
\usepackage[pagebackref]{hyperref}
\usepackage{pdfsync}
\usepackage{fancyhdr}
\usepackage{booktabs}
\usepackage{appendix,color}

\newcommand{\Circ}{\hspace{-0.5pt}\begin{minipage}[t]{9pt}\vspace{-6.5pt} \mbox{\fontsize{18}{18}\selectfont $\circ$}\end{minipage}\hspace{-1pt}}
\newcommand{\Bullet}{\hspace{-0.5pt}\begin{minipage}[t]{9pt}\vspace{-6.5pt} \mbox{\fontsize{18}{18}\selectfont $\bullet$}\end{minipage}\hspace{-0.5pt}}
\newcommand{\Square}{\begin{minipage}[t]{8.5pt}\vspace{-6.0pt} \mbox{\fontsize{9.5}{9.5}\selectfont $\square$}\end{minipage}\hspace{-0.5pt}}
\newcommand{\Blacksquare}{\begin{minipage}[t]{8.5pt}\vspace{-6.0pt} \mbox{\fontsize{9.5}{9.5}\selectfont $\blacksquare$}\end{minipage}\hspace{-0.5pt}}

\newcommand{\openbox}{\leavevmode
  \hbox to.77778em{%
  \hfil\vrule
  \vbox to.675em{\hrule width.6em\vfil\hrule}%
  \vrule\hfil}}

\newcommand{\bigo}{\mathcal O}

\newcommand{\lsdef}{\stackrel{def}{=}}

\newcommand{\cE}{\mathcal{E}}
\newcommand{\cO}{{\mathcal{O}}}
\newcommand{\cU}{\mathcal{U}}

\newcommand{\cL}{\mathcal{L}}
\newcommand{\cR}{\mathcal{R}}
\newcommand{\cG}{\mathcal{G}}
\newcommand{\cK}{\mathcal{K}}

\newcommand{\Ss}{\mathsf{s}}
\newcommand{\Sb}{\mathsf{b}}
\newcommand{\sB}{\mathsf{B}}
\newcommand{\sI}{\mathsf{I}}
\newcommand{\sS}{\mathsf{S}}
\newcommand{\sT}{\mathsf{T}}
\newcommand{\sC}{\mathsf{C}}
\newcommand{\sD}{\mathsf{D}}
\newcommand{\sM}{\mathsf{M}}

\newcommand{\bA}{\mathbf{A}}

\newcommand{\bD}{\mathbf{D}}

\newcommand{\bG}{\mathbf{G}}
\newcommand{\bH}{\mathbf{H}}

\newcommand{\bL}{\mathbf{L}}

\newcommand{\bR}{\mathbf{R}}
\newcommand{\bS}{\mathbf{S}}
\newcommand{\bU}{\mathbf{U}}

\newcommand{\bSigma}{\mathbf{\Sigma}}

\renewcommand{\b}{\boldsymbol}

\newcommand{\blockm}[4]{ \begin{pmatrix} #1   &{#2}  \\ #3   &{#4}  \end{pmatrix}}
\newlength{\numberwidth}
\newlength{\circlewidth}
\newcommand{\circlednumber}[1]{%
    \settowidth{\numberwidth}{\footnotesize$#1$}%
    {\large$\bigcirc$}%
    \hspace{-0.5\circlewidth}\hspace{-0.5\numberwidth}%
    {\raisebox{0.5pt}{\footnotesize$#1$}}%
    \hspace{0.5\circlewidth}\hspace{-0.5\numberwidth}%
}

\newenvironment{lsitemizei}{
    \begin{list}{$\bullet$}{
        \setlength{\itemsep}{0pt}
        \setlength{\parsep}{2pt}
        \setlength{\topsep}{2pt}
    }
}{\end{list}}

\newenvironment{lsitemizeii}{
    \begin{list}{$\diamond$}{
        \setlength{\itemsep}{0pt}
        \setlength{\topsep}{0pt}
    }
}{\end{list}}

\newtheorem{thm}{Theorem}
\newtheorem{property}{Property}
\newtheorem{cor}{Corollary}
\newcommand{\proof}[2][Proof]{\noindent {\it \textbf{#1}:} \ #2\ \hspace*{\fill} $\square$}

\newcounter{songquestion}

\newcounter{ericquestion}

\newcounter{generalnote}

\begin{document}

\title{Extension and optimization of the FIND algorithm: computing Green's and less-than Green's functions (with technical appendix)}

\author{S.~Li$^a$, E.~Darve$^{a,b}$ \\
$^a$Institute for Computational and Mathematical Engineering, Stanford University, \\
496 Lomita Mall, Durand building, Stanford, CA 94305, USA, \\
$^b$Department of Mechanical Engineering, Stanford University, \\
496 Lomita Mall, Durand Building, Room 209, Stanford, CA 94305, USA, \\
darve@stanford.edu}

\maketitle

{\bf Abstract.} \quad
The FIND algorithm is a fast algorithm designed to calculate certain entries of the inverse of a sparse matrix. Such calculation
is critical in many applications, e.g., quantum transport in nano-devices. We extended the algorithm to other matrix inverse
related calculations. Those are required  for example to calculate the less-than Green's function and the current density through
the device. For a 2D device discretized as an $N_x\times N_y$ mesh, the best known algorithms have a running time of
$\bigo(N_x^3N_y)$, whereas FIND only requires $\bigo(N_x^2N_y)$. Even though this complexity has been reduced by an order of
magnitude, the matrix inverse calculation is still the most time consuming part in the simulation of transport problems. We could
not reduce the order of complexity, but we were able to significantly reduce the constant factor involved in the computation
cost. By exploiting the sparsity and symmetry, the size of the problem beyond which FIND is faster than other methods typically
decreases from a $130\times130$ 2D mesh down to a $40\times40$ mesh. These improvements make the optimized FIND algorithm even
more competitive for real-life applications.

\newpage
\noindent{\large \textbf{Main Notations}}

\medskip

\noindent\begin{tabular}{lp{280pt}}
  \toprule
    $\sM$         & the set of all the nodes in the mesh, p.~\pageref{def:sM}  \\
    $\sT$         & the mesh is subdivided using nested dissection; $\sT$ is the resulting tree of clusters,
                    p.~\pageref{def:sT}     \\
    $\sC$, $\sC_g$ & $\sC$ is a set of mesh nodes after the partitioning. It is sometimes called a cluster.
                            It may designate a cluster at any level in the tree. A specific cluster is denoted as $\sC_g$, where $g$ is the
                            label used when partitioning the whole mesh, as shown in Fig.~\ref{fig:partitions-positive}. \\
	$\bar{\sC}$, $\sC_{-g}$ & $\bar{\sC}$ is the complement of $\sC$  with respect to $\sM$ and
                            $\sC_{-g} = \bar{\sC}_g$ as shown in Fig.~\ref{fig:partitions-negative}. p.~\pageref{def:sC} \\
	$\sB$, $\sB_g$ & set of boundary nodes of a cluster, $\sB_g$ is the set boundary nodes of $\sC_g$ and sometimes called the
                        boundary set of $\sC_g$. p.~\pageref{def:sB} \\
    $\sB_L$, $\sB_R$    & sets of boundary nodes originated from the left child cluster and the right cluster, p.~\pageref{def:BL-BR}  \\
    $\sI$, $\sI_g$        & set of inner nodes. $\sI_g$ is the set of inner nodes of $\sC_g$.
                            $\sI_g\cap\sB_g = \emptyset$ and $\sI_g\cup\sB_g = \sC_g$. p.~\pageref{def:sI}                \\
    $\sS$, $\sS_g$ & set of private inner nodes. $\sS_g$ is the set of private inner nodes of cluster $\sC_g$.
                        $\sS_g \subset \sI_g$.
                        p.~\pageref{def:sS}   \\
    $\sS_L$, $\sS_R$    & private inner nodes originated from the left child cluster and right cluster, p.~\pageref{def:SL-SR}  \\
    $\Ss$, $\Sb$       & size of the set $\sS$ and $\sB$, p.~\pageref{def:size}             \\
    $N, N_x, N_y$  & size of the mesh from the discretization of 2D device, p.~\pageref{def:Nx-Ny}.
                For square mesh, $N_x=N_y$ and we write it as $N$, p.~\pageref{def:N}.   \\
    $n$         & size of square matrices $\bA$, $\bSigma$, $\bG$, and $\bG^<$. For an $N_x\times N_y$ mesh, $n=N_xN_y$. p.~\pageref{def:n}  \\
    $\bA$, $\bA^{\dagger}$       & the sparse matrix from the discretization of the device and its conjugate transpose, p.~\pageref{def:A}  \\
    $\bG^r, \bG^<$     & matrix associated with the retarded Green's function
				p.~\pageref{def:Gr}\\
    $\bSigma$   & matrix associated with the self energy, p.~\pageref{def:Sigma}  \\
%
    $\bL$, $\bU$  & the LU factors of the sparse matrix $\bA$, p.~\pageref{def:LU}  \\
	$\bL_g$ & we decompose the LU factorization of $\bA$ into multiple partial factorizations: $\bL = \prod_g\bL_g$.
p.~\pageref{def:bLg}   \\
	$\bA_{g}$, $\bA_{g+}$ & the matrices in the elimination process right before and after eliminating $\sS_g$, respectively.
                            p.~\pageref{def:Ag-Agplus} \\
	$\bSigma_g$, $\bSigma_{g+}$ & we keep updating $\bSigma$ as we factorize $\bA$ and use a subscript to indicate
                            the stage in the process.  \\
    $\cL$, $\cL_g$  & $\cL_g = \bL_g(\sB_g, \sS_g)$ is the only non-zero off-diagonal block of $\bL_g$ used for the update step for
                        eliminating $\sS_g$.
                        Sometime it is simplified as $\cL$ when the subscript $g$ is not important. p.~\pageref{def:Lg} \\
    $\cU$, $\cU_g$       & $\cU_g$ is a simplified notation for $\bA_{g+}(\sB_g, \sB_g)$. We keep computing $\cU_g$ when we
                            calculate $\bG^r$.
                            Sometimes $\cU_g$ is further simplified as $\cU$ when the subscript $g$ is not important. p.~\pageref{def:Ug}     \\
    $\cR$, $\cR_g$       & $\cR_g$ is a simplified notation for $\bSigma_{g+}(\sB_g, \sB_g)$. We keep computing $\bR_g$ when we
                            calculate $\bG^<$.
                            Sometimes $\cR_g$ is further simplified as $\cR$ when the subscript $g$ is not important. p.~\pageref{def:Rg} \\
    \bottomrule
  \end{tabular}
\newpage

\section{Introduction}

In recent years, the non equilibrium Green's function formalism~\cite{datta97,anant02,ghosh2004molecular,martinez2007self}
(NEGF) has been widely used to study the properties of nanoscale MOS (metal-oxide-semiconductor) transistors as well as nanowires
and molecular electronic devices in which quantum effects are important. This formalism leads to the definition of two types of
Green's function, the retarded Green's function and less-than Green's function. Their discretized approximations satisfy the
following equations~\cite{li08}:
\begin{gather*}
	\bA = \bH_0(\cE)+\bSigma^r(\cE) \label{def:A}    \\
	\bG^r(\cE) = \bA^{-1} \label{def:Gr}   \\
	\bG^<(\cE) = \bG^r(\cE) \; \bSigma(\cE) \; [\bG^r]^\dagger(\cE)    \label{def:Gless}
\end{gather*}
where $\cE$ is the energy\label{def:E}, $\bH_0$ is the Hamiltonian representing a single particle, $\bSigma^r$ \label{def:Sigma}
is the self-energy used to model the contact of the nano-device (nano-transistors) with the (infinite) source and drain, and
$[\bG^r]^\dagger$ is the conjugate transpose of $\bG^r$. $\bSigma$ is the scattering self-energy matrix. It is typically assumed that
$\bSigma$ is a diagonal matrix, even though we will see later on that this assumption can be relaxed for our approach. The charge
density can be obtained from the diagonal entries of $\bG^<(\cE)$. In a typical calculation, one needs to solve a coupled
Schr\"{o}dinger-Poisson system of equations using some iterative scheme. Each calculation of the diagonal of $\bG^<(\cE)$ is
relatively fast, however the total computational cost can be very large since the charge density needs to be computed at all
energy levels $\cE$, and at each iteration of the Schr\"{o}dinger-Poisson equation solver.

In this paper, we will focus on discretization schemes for the Hamiltonian operator $\bH_0$ that use a 2D Cartesian grid. The
method can be extended to more general schemes including finite-element but the description is easier in the context of Cartesian
grids. We will further assume that a 5 point finite-difference stencil is used, although the method can be extended to arbitrary
stencils. This stencil is the most common for these types of problems. Finally even though the matrix is not symmetric we require
that $\bA_{ij} \neq 0 \; \Leftrightarrow \; \bA_{ji} \neq 0$. All these requirements can be removed but these assumptions will
simplify the discussion of the method.

In previous publications~\cite{darve04,li07,li08}, we discussed an efficient method, FIND (Fast Inverse using Nested
Dissection), to calculate the diagonal of $\bG^r=\bA^{-1}$. The novel contributions in this paper include:
\begin{itemize}
	\item Extension of FIND~\cite{darve04,li07,li08} for computing the diagonal of $\bG^< = \bA^{-1} \bSigma
\bA^{-\dagger}$. See Section~\ref{extension-diagonal}.
	\item Extension for computing the off-diagonal entries of $\bG^r$. This is required for computing the current density
(see for example~\cite{lake1997single}, eq.~(20) pp~7847). This is presented in Section~\ref{sec:current}.
	\item Optimizations using certain sparsity patterns in $\bA$ to further reduce the computational cost
(Section~\ref{sec:extra}).
	\item Optimizations using the symmetry of $\bA$ (Section~\ref{symmetry}).
	\item Optimizations specific to computing the off-diagonal entries for current density calculations.
	\item Proof of symmetry and/or positive definiteness for matrices arising in the algorithm, when $\bA$ is itself
symmetric and/or positive definite. See~\ref{symmetry}.
	\item Detailed analysis of the computational cost and storage requirement for the different approaches.
\end{itemize}
Pseudo-codes are given for the key algorithms. Numerical results with accuracy and performance benchmarks are given at the end (Section~\ref{sec:num}). The performance improvements vary but are on the order of 2x in general. In the appendix, we derive a series of technical results that are required to prove some of the statements in this paper.

Before beginning, we make some general comments regarding our approach and how it relates to other techniques. In~\cite{borm2003hierarchical,hackbusch99}, methods are proposed to calculate the inverse of a dense matrix, of size $N\times N$, in $O(N \ln^2 N)$ steps using low rank approximations of certain off-diagonal blocks. Similar techniques include ${\cal H}^2$-matrices, quasi-separable matrices, semi-separable matrices, and HSS-matrices~\cite{borm2003hierarchical,chandrasekaran06a,chandrasekaran06b,chandrasekaran08}. Mamaluy et al.'s method uses the special dependency of the problem on $\cE$ to reduce the computational cost~\cite{mamaluy03,khan2007quantum,sabathil2004prediction,mamaluy2004contact,mamaluy2004electron,birner-ballistic}. This last method requires computing the eigenvectors and eigenvalues of the decoupled Hamiltonian. In~\cite{tang2010probing}, a method is proposed in cases where the inverse matrix exhibits a decay property, i.e., when many of the entries of the inverse are small. In this paper, we seek to find exact methods that do not require special matrix properties or solving an eigenvalue problem.

An efficient approach in this second class of methods is the recursive Green's function (RGF) method~\cite{svizhenko02} (see also the earlier work of Erisman et al.~\cite{erisman1975computing}) with running time of order $\bigo(N_x^3N_y)$, where $N_x$ and $N_y$ \label{def:Nx-Ny} are the number of points on the grid in the $x$ and $y$ direction after discretization, respectively ($N_x < N_y$). This technique was generalized in~\cite{petersen08} to block tri-diagonal matrices. RGF type methods perform a forward recurrence, based on LU factorizations. This produces a single entry in the inverse. From this entry, a backward recurrence progressively computes the remaining entries. In contrast, FIND~\cite{darve04,li07,li08} performs two elimination passes, based on LU factorizations, typically described as upward elimination (similar to the LU factorization in RGF), and a downward elimination (which has no equivalent in RGF). At the end of the process, the inverse is directly computed as a by-product of the downward elimination. The work of Lin et al.~\cite{lin2009fast,lin2011selinv} follows an approach similar to FIND.

As a note, we mention that the algorithm presented below uses a decomposition of the problem based on nested dissection. The original nested dissection of George et al.~\cite{george73}, which is used to solve linear systems, requires defining {\bf separator sets.} A separator set is a set of mesh nodes such that it divides the mesh in two sets $\sC_1$ and $\sC_2$, and $\bA_{ij} = 0$ if $i \in \sC_1$ and $j \in \sC_2$. An important difference between nested dissection and FIND is that, in the original nested dissection, a separator set is ``shared'' by the neighboring clusters. In FIND, since we need more independence among clusters, each cluster has its own boundary set that separates it from other clusters. These boundary sets are similar to the separator sets (generally called width-1 separators) used in nested dissection, but form so-called width-2 separators, with basically one ``separator'' for each of the two neighboring clusters. See for example~\cite{liu1974solution,george1983row,gilbert1987parallel,pothen1990partitioning} for a definition of width-$l$ separators. Width-$l$ separators with $l>1$ are sometimes called wide separators. However, it is possible to construct a FIND algorithm that uses the more classical width-1 separators. Currently, FIND belongs to the class of {\bf right-looking} algorithms (see~\cite{heath1991parallel} page 426, 427). A variant of FIND requiring only width-1 separators can be created using a {\bf left-looking} strategy for the additions.\footnote{Specifically, this variant uses right-looking multiplications and left-looking additions.} Even though the latter might be more efficient, in this paper, we focus on the original FIND scheme, which uses a right-looking approach with width-2 separators.

\section{Previous work}

\subsection{Brief description of the FIND algorithm for $\bG^r = \bA^{-1}$}

\label{brief}

The basic idea of the FIND algorithm is to simultaneously perform many LU factorizations on the given matrix to obtain the
diagonal elements of its inverse. By performing the LU factorizations in a specific order that minimizes fill-ins
\cite{george73}, we preserve most of the sparsity of the given matrix $\bA$, and thus make the LU factorizations very efficient.
Once an LU factorization is complete, we can easily compute the last entry on the diagonal of the inverse: for an $n\times
n$\label{def:n} matrix $(\bG^r)^{-1} = \bL\bU$,\label{def:LU} we have $\bG^r_{nn} = 1/\bU_{nn}$. Although we can compute only the
$(n, n)$ entry of the diagonal by this formula, we can choose any node and reorder the original matrix to make that node
correspond to the $(n, n)$ entry of the reordered matrix. In this way, all the diagonal elements of $\bG^r$ can be computed.

If we have to perform a full LU factorization for each of the $n$ reordered matrices, the algorithm will not be computationally
efficient even though each LU factorization is very fast. However, we can decompose each LU factorization into partial LU
factorizations: $\bA=\prod_{g}\bL_g\prod_{g}\bU_g$, where $\bL_g$ corresponds to the elimination of $\sS_g$ and $\sS_g$ is a set
of private inner nodes that we will discuss in Section~\ref{sec:basic-detailed}.\label{def:bLg} Due to the sparsity pattern of
$\bA$, we can perform them independently and moreover, if we reorder $\bA$ properly, many partial factorizations for different
reordered matrices turn out to be identical. As a result, we can reuse the partial results multiple times for different orderings
thereby reducing considerably the computational cost.

\begin{figure}[htbp]
\centering

\subfigure[Partition into basic clusters.]{ \vspace{-15pt}
\begin{minipage}[b]{175pt}
\centering
\includegraphics[width=75pt]{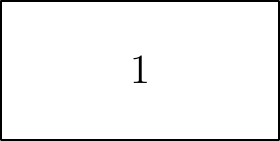} \hspace{2pt}
\includegraphics[width=75pt]{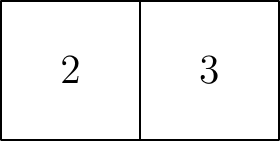} \\
\vspace{10pt}
\includegraphics[width=75pt]{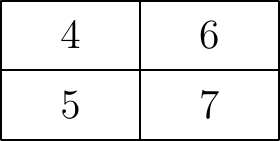} \hspace{2pt}
\includegraphics[width=75pt]{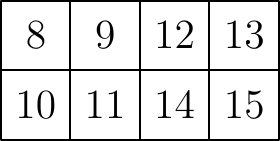}
\end{minipage}
\label{fig:partitions-positive}} \hspace{15pt} \subfigure[Partition into basic and complement clusters.]{
\vspace{-15pt}
\begin{minipage}[b]{100pt}
\centering
\includegraphics[width=75pt]{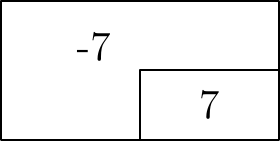} \\
\vspace{10pt}
\includegraphics[width=75pt]{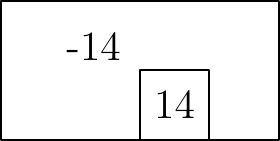}
\end{minipage}
\label{fig:partitions-negative}
}

\caption{Partitions of a 2D rectangular mesh into clusters} \label{fig:partitions}
\end{figure}

More specifically, we partition the mesh into subsets of mesh nodes. Fig.~\ref{fig:partitions-positive} shows the partitions of
the entire mesh into 2, 4, and 8 clusters, for a 2D rectangular grid. Then, these subsets are merged to form larger clusters.
Each merge corresponds to a partial Gaussian elimination. The clusters shown on the left panel~\ref{fig:partitions-positive} are
called \emph{basic clusters} and usually denoted as $\sC$.\label{def:sC} They are labeled with positive integers therefore
sometimes called \emph{positive clusters}. For reasons we will see soon, we also consider the complement of the basic clusters
with respect to the whole mesh $\sM$.\label{def:sM} These clusters are called \emph{complement clusters} and denoted as
$\bar{\sC}$.\label{def:sC-bar} They are labeled with negative integers therefore sometimes called \emph{negative clusters} as
shown in Fig.~\ref{fig:partitions-negative}.

\begin{figure}[htbp]
\centering

\subfigure[The basic cluster tree of a mesh.]{ \vspace{-20pt}
\begin{minipage}[b]{163pt}
\centering
	\includegraphics[width=162pt]{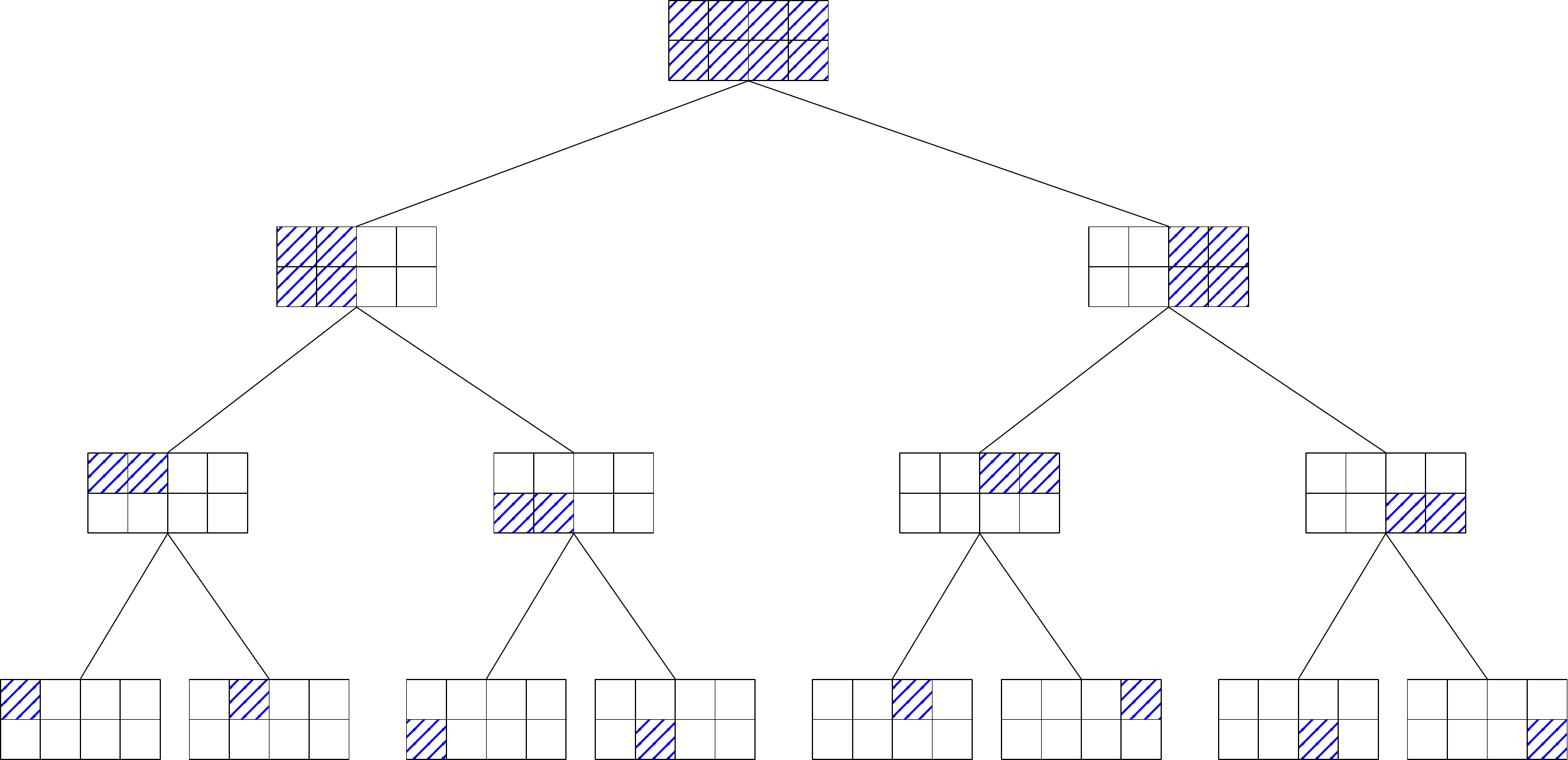}
	\label{fig:cluster-tree-upward}
\end{minipage}
} \hspace{-1pt} \subfigure[The complement-cluster tree of a mesh.]{ \vspace{-20pt}
\begin{minipage}[b]{163pt}
	\includegraphics[width=162pt]{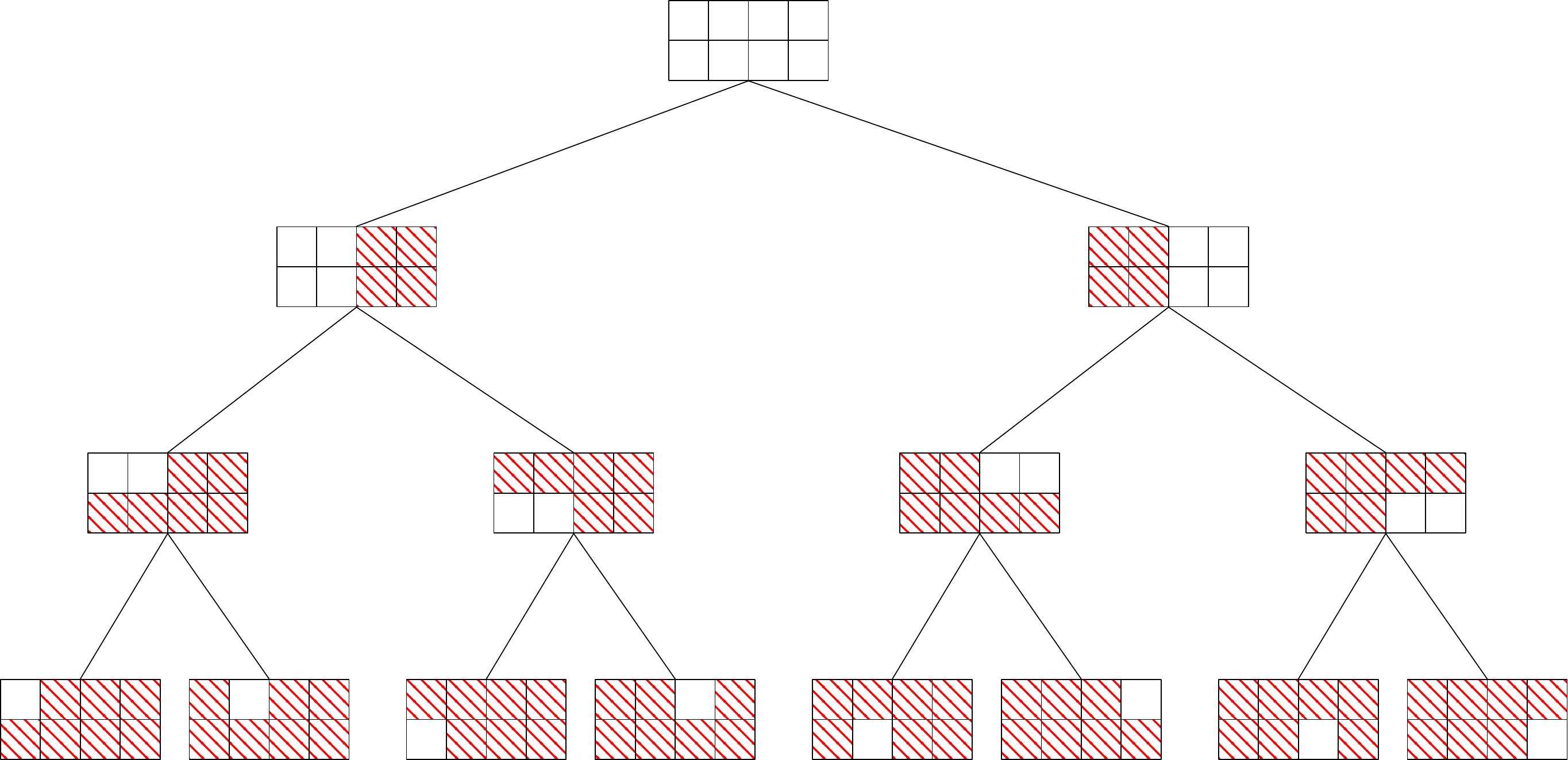}
	\label{cluster-tree-downward}
\end{minipage}
}
\caption{Cluster trees.}
\end{figure}

We organize these basic and complement clusters in trees to see how the partial eliminations are related to the full elimination.
Fig.~\ref{fig:cluster-tree-upward} shows how the basic clusters form a \emph{basic cluster tree}, usually denoted as
$\sT$.\label{def:sT} For each of the leaf clusters (\emph{target cluster}) in the figure, we order the matrix $\bA$, such that
the entries corresponding to the target cluster (\emph{target entries}) form a block at the end of the matrix. We then eliminate
all the columns corresponding to the complement of the target cluster to compute the target entries in $\bA^{-1}$.

Because the complements of the target clusters overlap significantly, we organize these clusters, together with other complement
clusters, into a similar \emph{complement cluster tree}, as shown in Fig.~\ref{cluster-tree-downward}. Unlike the clusters in
Fig.~\ref{fig:cluster-tree-upward}, a parent cluster in Fig.~\ref{cluster-tree-downward} is not a union of its two child
clusters, but rather their intersection.

\begin{figure}[htbp]
\begin{minipage}[b]{170pt}
\includegraphics[width=180pt]{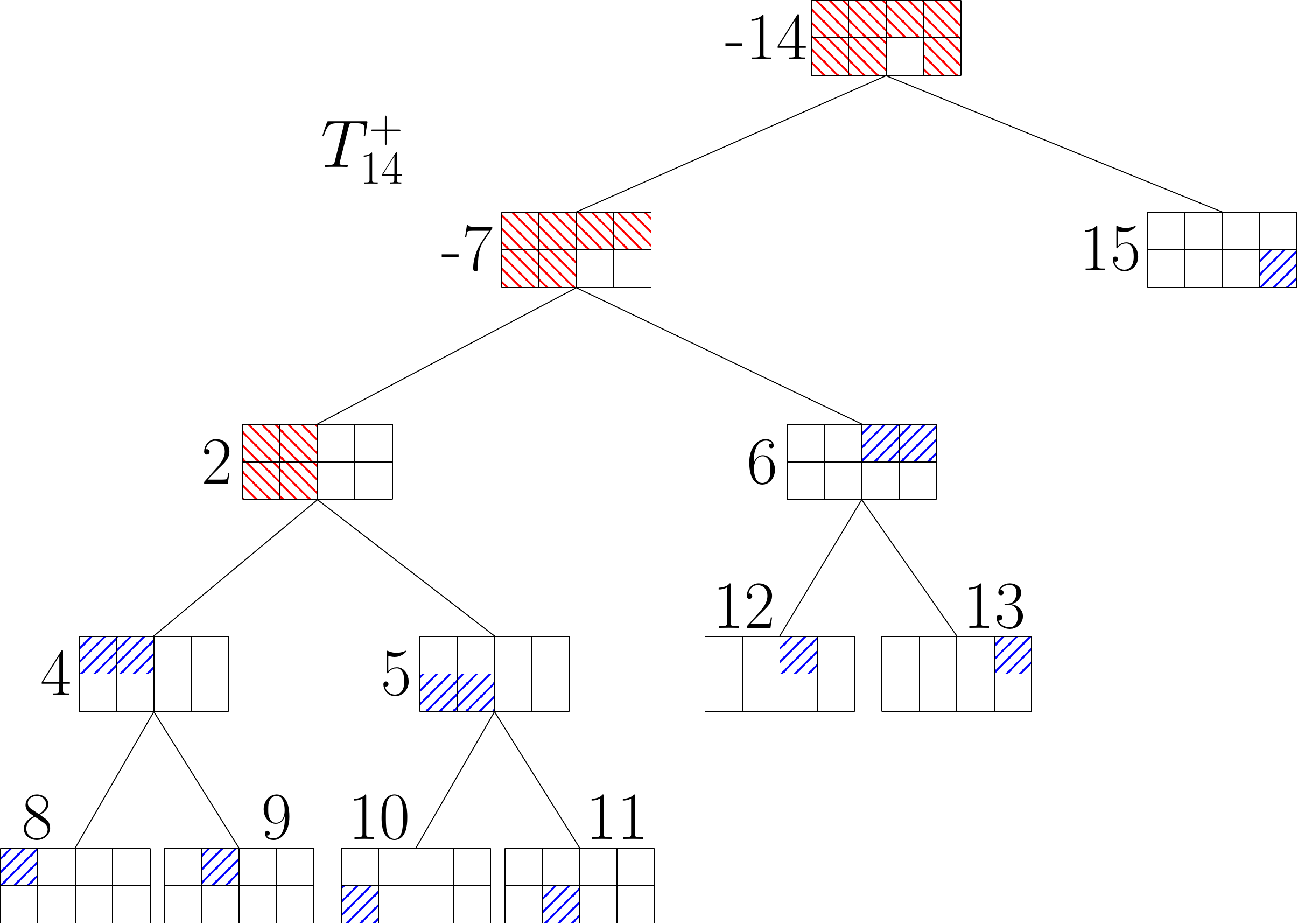}
\end{minipage} \nolinebreak
\begin{minipage}[b]{170pt}
\includegraphics[width=180pt]{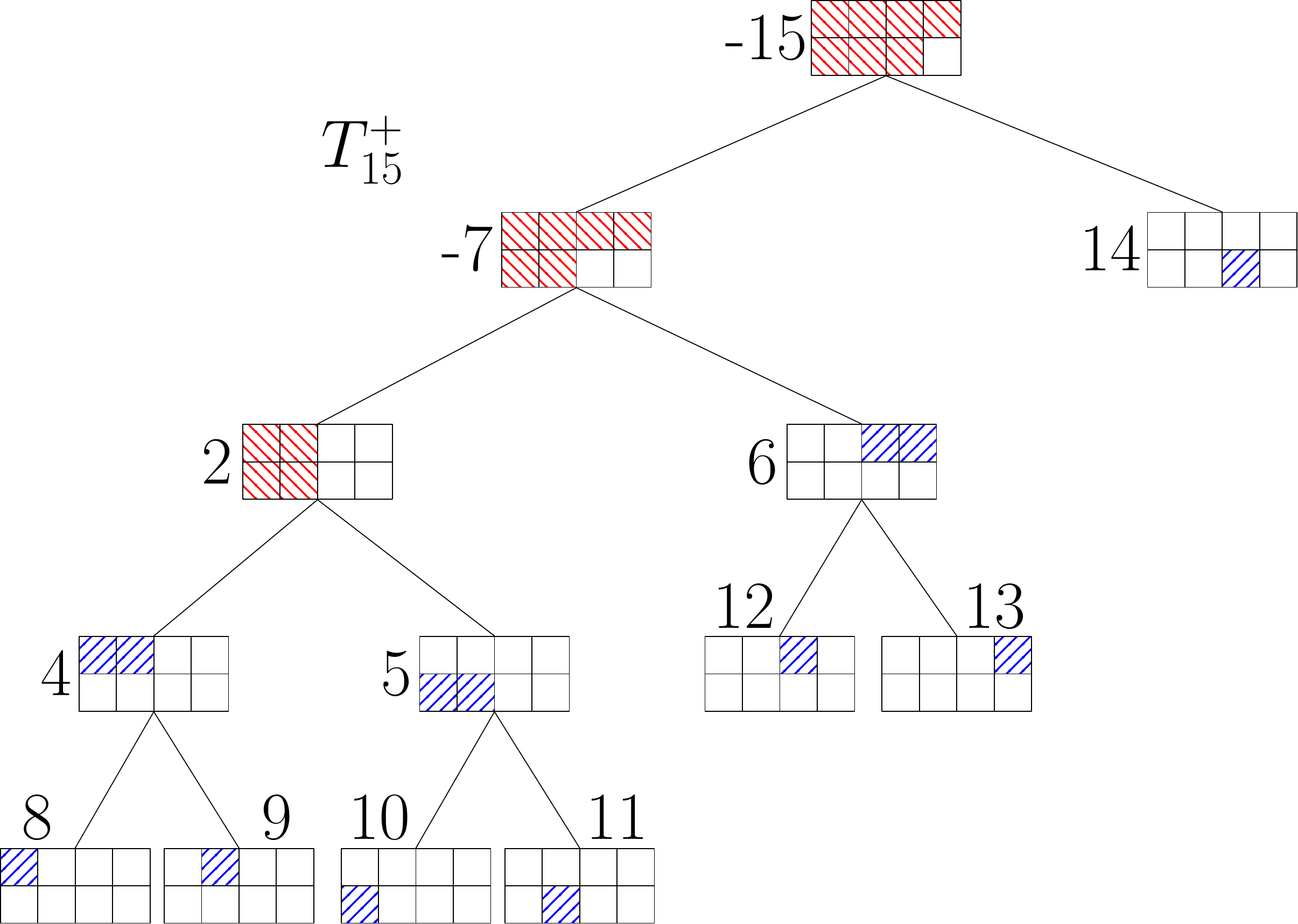}
\end{minipage}
\caption{ Two augmented trees with target clusters $14$ and $15$.}
\label{two-augmented-tree}
\end{figure}

The two cluster trees in Fig.~\ref{two-augmented-tree} show how two leaf complement clusters $-14$ and $-15$ are computed by the
FIND algorithm. Such trees are called \emph{augmented trees} and are denoted as $\sT_r^+$ with $r$ indicating the target cluster,
i.e., the complement of the root cluster. Each tree is associated with an ordering of $\bA$ with the target entries appearing at
the end of the matrix. The path $(-3)$-$(-7)$-$(-14)$ in $\sT_{14}^+$ and the path $(-3)$-$(-7)$-$(-15)$ in $\sT_{15}^+$ are part
of the complement-cluster tree; the other subtrees are copies from the basic cluster tree. We can see that these two trees are
almost identical, which means that the two full LU factorizations of $\bA$ with different orderings significantly overlap.

\begin{figure}[htbp]
\flushright
\includegraphics[width=280pt]{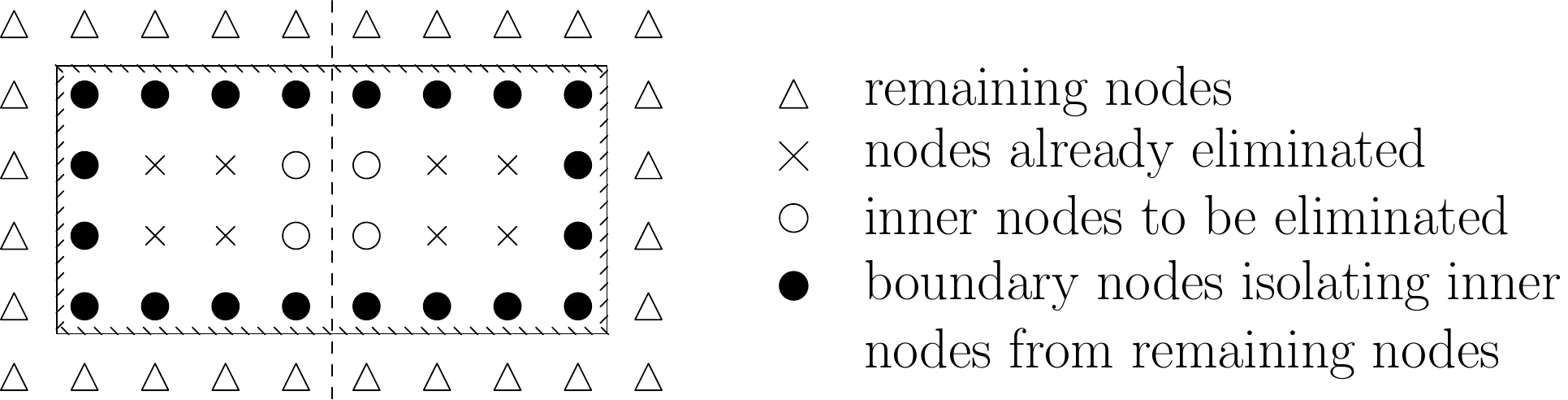}
\caption{The nodes within the rectangular frame form a large cluster. Two child clusters are shown on the left and right.}
\label{fig:simple-merge}
\end{figure}

At each level of the two trees, two clusters are merged into a larger one. We define the \emph{inner nodes} of a cluster $\sC$ as
the subset of nodes $\sI$ in a cluster such that $\bA_{ij}=0$ for any $i \in \sI$ and $j \notin \sC$. Each merge is equivalent to
the Gaussian elimination of all the inner nodes of the resulting cluster. For example, in Fig.~\ref{fig:simple-merge}, the large
cluster has been surrounded by a solid line with tics. The two children clusters are on the left and right. During earlier steps
in the method, the $\times$ nodes have been eliminated. At the current step, the $\circ$ nodes are being eliminated. As a result,
all the inner nodes of the clusters have been eliminated at the end of this operation.

Such a merge corresponds to an independent partial LU factorization since there is no connection between the eliminated nodes
($\circ$) and the remaining nodes ($\vartriangle$). This is better seen below in the Gaussian elimination in~\eqref{update-rule}
for the merge in Fig.~\ref{fig:simple-merge}:\\
\begin{minipage}[b]{100pt}
\begin{equation*}
\begin{pmatrix}
\bA(\circ, \circ) &\hspace{-6pt} \bA(\circ, \bullet)  &\hspace{-6pt} \b 0  \\
\bA(\bullet, \circ) &\hspace{-6pt} \bA(\bullet, \bullet)  &\hspace{-6pt} \bA(\bullet,
\vartriangle)   \\
\b 0   &\hspace{-6pt} \bA(\vartriangle, \bullet)   &\hspace{-6pt} \bA(\vartriangle, \vartriangle)
\end{pmatrix}
\Rightarrow
\begin{pmatrix}
\bA(\circ, \circ) &\hspace{-6pt} \bA(\circ, \bullet)  &\hspace{-6pt} \b 0  \\
\b 0 &\hspace{-6pt} \bA^*(\bullet, \bullet)  &\hspace{-6pt} \bA(\bullet,\vartriangle)   \\
\b 0   &\hspace{-6pt} \bA(\vartriangle, \bullet)   &\hspace{-6pt} \bA(\vartriangle, \vartriangle)
\end{pmatrix},
\end{equation*}
\end{minipage}
\begin{equation}\hspace{-0pt} \mbox{where }\bA^*(\bullet, \bullet) = \bA(\bullet, \bullet) - \bA(\bullet, \circ)\bA(\circ,
\circ)^{-1}\bA(\circ, \bullet)\label{update-rule}
\end{equation}
\normalsize In~\eqref{update-rule}, $\circ$ and $\bullet$ each corresponds to all the $\circ$ nodes and $\bullet$ nodes in
Fig.~\ref{fig:simple-merge}, respectively. $\bA^*(\bullet, \bullet)$ is the result of the partial LU factorization. It will be
stored for reuse when we go through the cluster trees.

To make the most of the overlap, we do not perform the eliminations for each augmented tree individually. Instead, we perform the
elimination based on the basic cluster tree and the complement-cluster tree in the FIND algorithm. We start by eliminating all
the inner nodes of the leaf clusters in Fig.~\ref{fig:cluster-tree-upward} followed by eliminating the inner nodes of their
parents recursively until we reach the root cluster. This is called the \emph{upward pass}. Once we have the partial elimination
results from upward pass, we can use them to traverse the complement-cluster tree level by level, from the root to the leaf
clusters. This is called the \emph{downward pass}.

\subsection{Formal description of the FIND algorithm} \label{sec:basic-detailed}

Section~\ref{brief} briefly illustrated how the algorithm works by showing that the eliminations of inner nodes are independent
of each other (Eq.~\eqref{update-rule}) and such eliminations are mostly shared among different augmented trees
(Fig.~\ref{two-augmented-tree}). Below we will show in a more formal way why the basic FIND algorithm works.

In Fig.~\ref{fig:simple-merge}, the $\bullet$ nodes are called the boundary nodes, denoted as $\sB_i$.\label{def:sB} The inner
nodes consist of the $\circ$ and $\times$ nodes, denoted as $\sI_i$.\label{def:sI} We call the $\circ$ nodes \emph{private inner
nodes} and denote them as $\sS_i$.\label{def:sS} Private inner nodes are the nodes that actually get eliminated when we process
the cluster. To guarantee the independence of the eliminations shown by \eqref{update-rule} and the overlap of the eliminations
for different target nodes illustrated by Fig.~\ref{two-augmented-tree}, we need to order the original matrix $\bA$ such that all
the columns corresponding to the private nodes of each cluster in Fig.~\ref{fig:partitions} are grouped together. It can be shown
that for each target cluster, the whole mesh can be partitioned as a series of sets of private inner nodes ($\sS_i$), the set of
boundary nodes $\sB_{\bar{r}}$ of the complement of the target cluster, and the target cluster $\sC_r$. We order the columns in
$\bA$ based on the ordering of these sets: $\sS_i$ of the leaf clusters first, followed by $\sS_i$ of their parents, $\sS_i$ of
the clusters in the next level in the augmented tree, $\ldots$, $\sS_{\bar{r}}$, $\sB_{\bar{r}}$, and finally $\sC_r$. We call
such ordering a \emph{consistent ordering} with respect to a target cluster. Note that the consistent ordering is not unique
since we do not require any specific ordering of clusters at the same level in the augmented tree. Here is an example of a
consistent ordering with respect the target cluster 14: $\sS_8$, $\sS_9$, $\sS_{10}$, $\sS_{11}$, $\sS_4$, $\sS_5$, $\sS_{12}$,
$\sS_{13}$, $\sS_{\bar{3}}$, $\sS_6$, $\sS_{\bar{7}}$, $\sS_{15}$, $\sS_{\bar{14}}$, $\sB_{\bar{14}}$, and $\sC_{14}$.

For notation convenience, we write $\sS_i<\sS_j$ if $\sS_i$ appears before $\sS_j$ in a ordering and write
$\cup_{\sS_j<\sS_i}\sS_j$ as $\sS_{<i}$ and $(\cup_{\sS_i<\sS_j}\sS_j) \cup \sB_{-r} \cup \sC_r$ as $\sS_{>i}$. We also denote as
$\bA_g$ the matrix $\bA$ after all the columns corresponding to the private nodes appearing before $\sS_g$ have been eliminated
and denote as $\bA_{g+}$ the matrix with all the columns before and including $\sS_g$ have been eliminated.\label{def:Ag-Agplus}
When necessary, we write them as $\bA_{r, g}$ and $\bA_{r, g+}$ to indicate explicitly the target cluster $r$ of the elimination
process. Now, we can rewrite the one-step elimination illustrated by~\eqref{update-rule} more formally as:
\begin{align*}
  \bA_g & =
    \begin{pmatrix}
      \bA_g(\sS_{<g}, \sS_{<g}) &\bA_g(\sS_{<g}, \sS_g)   &\bA_g(\sS_{<g}, \sB_g)   &\bA_g(\sS_{<g}, \sS_{>g}\backslash \sB_g)  \\
      \b0                   &\bA_g(\sS_g, \sS_g)      &\bA_g(\sS_g, \sB_g)      &\b0  \\
      \b0                   &\bA_g(\sB_g, \sS_g)      &\bA_g(\sB_g, \sB_g)   &\bA_g(\sB_g, \sS_{>g}\backslash \sB_g)   \\
      \b0                   &\b0                  &\bA_g(S_{>g}\backslash \sB_g, \sB_g) &\bA_g(\sS_{>g}\backslash \sB_g, S_{>g}\backslash \sB_g)
    \end{pmatrix}
\intertext{$\b \Rightarrow$}
  \bA_{g+} & =
    \begin{pmatrix}
      \bA_g(\sS_{<g}, \sS_{<g}) &\bA_g(\sS_{<g}, \sS_g)   &\bA_g(\sS_{<g}, \sB_g)  &\bA_g(\sS_{<g}, \sS_{>g}\backslash \sB_g)  \\
      \b0   &\bA_g(\sS_g, \sS_g)  &\bA_g(\sS_g, \sB_g) &\b0  \\
      \b0   &\b0   &\bA_{g+}(\sB_g, \sB_g) &\bA_g(\sB_g, \sS_{>g}\backslash \sB_g)   \\
      \b0   &\b0  &\bA_g(\sS_{>g}\backslash \sB_g, \sB_g) &\bA_g(\sS_{>g}\backslash \sB_g, \sS_{>g}\backslash \sB_g)
    \end{pmatrix}
\label{theorem1:matrices}
\end{align*}
where
\begin{equation}
\bA_{g+}(\sB_g, \sB_g) = \bA_g(\sB_g, \sB_g) - \bA_g(\sB_g, \sS_g)\bA_g(\sS_g, \sS_g)^{-1}\bA_g(\sS_g,
\sB_g).\label{eq:major-update}
\end{equation}
For notation simplicity, we write $\bA_{g+}(\sB_g, \sB_g)$ as $\cU_g$.\label{def:Ug} This is the block we keep computing in our
algorithm. Also recall that $\bA=\prod_{g}\bL_g\prod_{g}\bU_g$, where $\bL_g$ corresponds to the elimination of $\sS_g$. So we
have $\bA_{g+} = \bL_g^{-1}\bA_g$ with $\bL_g^{-1}(\sB_g, \sS_g) = -\bL_g(\sB_g, \sS_g) = -\bA_g(\sB_g, \sS_g)\bA_g(\sS_g,
\sS_g)^{-1}$ the only non-zero off-diagonal block of $\bL_g$ (see Li et al. for more details)~\cite{li08}. We can further define
\begin{equation}
	\cL_g = \bL_g(\sB_g, \sS_g),\label{def:Lg}
\end{equation}
and then we can write~\eqref{eq:major-update} as
\begin{equation}
\bA_{g+}(\sB_g, \sB_g) = \bA_g(\sB_g, \sB_g) - \cL_g \bA_g(\sS_g,
\sB_g).
\end{equation}

Note that in the above matrices, for notation convenience, $\bA_g(\bullet, \sB_g)$ \label{block-notation} is written as a block.
In reality, however, it is usually not a block in $\bA$ of any ordering in our algorithm because the columns of $\bA$ are ordered
based on $\sS_g$ but not $\sB_g$. It can be shown that the matrix preserves the sparsity pattern required by the above formula
during the elimination process, which leads to the following theorem:
\begin{thm} \label{thm:FIND-major}
For any target clusters $r$ and $s$ such that $\sC_g\in \sT_{r}^+$ and $\sC_g\in \sT_{s}^+$, we have
\begin{equation*}
  \bA_{r, g+}(\sB_g, \sB_g) = \bA_{s, g+}(\sB_g, \sB_g).
\end{equation*}
\end{thm}
Theorem~\ref{thm:FIND-major} shows that the partial elimination results are common for matrices with different orderings during
the elimination process. This is the key foundation of the basic FIND algorithm. For more details, please see \cite{li08}.

\subsection{Cost analysis of the FIND algorithm for $\bA^{-1}$}

For analytic simplicity, we assume that $\bA$ comes from a square grid of size $N\times N$ \label{def:N}.
In this case, the total cost is $\bigo(N^3)$. This is the same computational cost in the $\bigo$ sense as the nested dissection algorithm of George et al.~\cite{george73}. It is now apparent that FIND has the same order of computational complexity as a single LU factorization, even though it calculates the diagonal of the inverse matrix.

Similar analysis tells us that the total memory cost is $\bigo(N^2\log(N))$. Both costs are asymptotically better than those for the best known algorithm, RGF, given by Svizhenko et al.~\cite{svizhenko02}. However, the constant factor in RGF~\cite{petersen09} is smaller than ours \cite{li08}. As a result, with the chosen tree decomposition, we are slower than RGF for meshes smaller than $130\times 130$. The single-layer-separator FIND variant discussed in the introduction would make FIND as fast as RGF even for small matrices. The following two sections will explain how to exploit the sparsity and symmetry to reduce the constant factor in FIND.

\section{Extension of the FIND algorithm}\label{sec:extension}

In addition to computing the diagonal entries of the matrix inverse, the algorithm can be extended to computing the diagonal
entries and certain off-diagonal entries of $\bG^< = \bA^{-1}\bSigma \bA^{-\dagger}$, which are required for the charge and
current densities.

\subsection{Computing the diagonal entries of $\bG^<$}\label{sec:compute-gless}

\label{extension-diagonal}

Intuitively, $\bA^{-1}\bSigma \bA^{-\dagger}$ can be computed in the same way as $\bA^{-1}$ because we have
\begin{equation}
		\extracolsep{-2em}
	\begin{array}{rr@{\,}l}
                  & \bG^< & = \bA^{-1}\bSigma \bA^{-\dagger}    \\
\Rightarrow \quad & \bA\bG^<\bA^\dagger &= \bSigma          \\
\Rightarrow \quad & \bU\bG^<\bU^\dagger &=  \bL^{-1}\bSigma \bL^{-\dagger} \\
\Rightarrow \quad & [\bG^<]_{nn} &=  (\bU_{nn})^{-1} (\bL^{-1}\bSigma \bL^{-\dagger})_{nn} (\bU_{nn})^{-\dagger}. \label{eq:gless-nn}
	\end{array}
\end{equation}

Now, it remains to show how to compute $(\bL^{-1}\bSigma \bL^{-\dagger})_{nn}$ efficiently. In the appendix, we show that if we order $\bSigma$ in a same way as $\bA$, the pattern of $\bSigma$ will be similar to
that of $\bA$ during the update process. Therefore, we can decompose the computation for the last block of $\bL^{-1}\bSigma
\bL^{-\dagger}$ into multiple steps with each step involving only a few small blocks. More specifically, we have a sequence of
matrices $\bSigma_g$ that start from $\bSigma$ and end at $(\bL^{-1}\bSigma \bL^{-\dagger})_{nn}$. This sequence is synchronized
with the updates on $\bA$, and similar to the update $\bA_{g+} = \bL_g^{-1}\bA_g$, we have $\bSigma_{g+} = \bL^{-1}_g \bSigma_g
\bL^{-\dagger}_g$. This one-step update is given by:
\begin{align}
\bSigma_g & =
\left(
{\footnotesize
\begin{array}{@{\,}cccc@{\,}}
\bSigma_g(\sS_{<g}, \sS_{<g}) &\bSigma_g(\sS_{<g}, \sS_g)   &\underline{\bSigma_g(\sS_{<g}, \sB_g)}   &\bSigma_g(\sS_{
<g}, \sS_{>g}\backslash \sB_g)  \\
\bSigma_g(\sS_g, \sS_{<g})    &\bSigma_g(\sS_g, \sS_g)      &\underline{\bSigma_g(\sS_g, \sB_g)}      &\b 0
                            \\
\underline{\bSigma_g(\sB_g, \sS_{<g})}    &\underline{\bSigma_g(\sB_g, \sS_g)}      &\underline{\bSigma_
g(\sB_g, \sB_g)}   &\bSigma_g(\sB_g, \sS_{>g}\backslash \sB_g)        \\
\bSigma_g(\sS_{>g}\backslash \sB_g)                   &\b 0     &\bSigma_g(\sS_{>g}\backslash
\sB_g, \sB_g) &\bSigma_g(\sS_{>g}\backslash \sB_g, \sS_{>g}\backslash \sB_g)
\end{array}
}\right)
\label{eq:Sigma-block-representation}
\end{align}
\begin{align*}
\intertext{$\Rightarrow$}
\bSigma_{g+} & = \left(
{\footnotesize
\begin{array}{@{\,}cccc@{\,}}
\bSigma_g(\sS_{<g}, \sS_{<g}) &\bSigma_g(\sS_{<g}, \sS_g)   &{\bSigma_{g+}(\sS_{<g}, \sB_g)}   &\bSigma_g(
\sS_{<g}, \sS_{>g}\backslash \sB_g)  \\
\bSigma_g(\sS_g, \sS_{<g})    &\bSigma_g(\sS_g, \sS_g)      &{\bSigma_{g+}(\sS_g, \sB_g)}      &\b 0  \\
{\bSigma_{g+}(\sB_g, \sS_{<g})}    &{\bSigma_{g+}(\sB_g, \sS_g)}      &{\bSigma_{g+}(\sB_g, \sB_g)}   &\bSigma_g(\sB_g, \sS_{>g}\backslash \sB_g)        \\
\bSigma_g(\sS_{>g}\backslash \sB_g)                   &\b 0 &\bSigma_g(\sS_{>g}\backslash
\sB_g, \sB_g) &\bSigma_g(\sS_{>g}\backslash \sB_g, \sS_{>g}\backslash \sB_g)
\end{array}
}\right),
\end{align*}
where
\begin{equation}
	\begin{split}
		\bSigma_{g+}(\sB_g, \sB_g) = &\ \bSigma_g(\sB_g, \sB_g) - \cL_g\bSigma_g(\sS_g, \sB_g) \\
		& - \bSigma_g(\sB_g, \sS_g)\cL_g^{\dagger} + \cL_g\bSigma_g(\sS_g, \sS_g)\cL_g^{\dagger},		
	\end{split}
\label{eq:gless-update-rule}
\end{equation}

The difference between the updates on $\bA$ and the updates on $\bSigma$ is that the latter will not only change $\bSigma_g(\sB_g, \sB_g)$, but also $\bSigma_g(\sS_{\leq g}, \sB_g)$ and $\bSigma_g(\sB_g, \sS_{\leq g})$. These blocks are underlined in \eqref{eq:Sigma-block-representation}. The reason for these changes is that we multiply $\bSigma$ by $\bL^{-1}$ from both sides and the subdiagonal blocks of $\bSigma_g$ are not cleared through the update process. The changes to $\bSigma_g(\sS_{\leq g}, \sB_g)$ and $\bSigma_g(\sB_g, \sS_{\leq g})$, however, do not invalidate the update process because in any future update step $h>g$, we are only interested in the operation on the submatrix $\bSigma_h(\sS_{\geq h}, \sS_{\geq h})$ and $\sS_{\leq g} \cap \sS_{\geq h} = \emptyset$. See the appendix for a rigorous description and proof.

Note that even though $\bSigma_{r, g+}$ generally depends on $r$ (corresponding to the ordering for different target clusters),
$\bSigma_{r, g+}(\sB_g, \sB_g)$ does not, as shown in the appendix. For notation
simplicity, we write $\bSigma_{g+}(\sB_g, \sB_g)$ as $\cR_g$.\label{def:Rg} The major part of the computation for $\bG^<$ is to
compute $\cR_g$ (for $\bG^<$), along with $\cU_g$ (for both $\bG^r$ and $\bG^<$), for all clusters.

\subsubsection{The algorithm and pseudo-codes}

Let $\sC_i$ and $\sC_j$ be the two children of $\sC_g$ in $\sT_r^+$. To compute $\cR_g$, we need both $\bA_g(\sS_g\cup \sB_g,
\sS_g\cup \sB_g)$ and $\bSigma_g(\sS_g\cup \sB_g, \sS_g\cup \sB_g)$. $\bA_g(\sS_g\cup \sB_g, \sS_g\cup \sB_g)$ follows the same
update rule as in the basic FIND algorithm; and $\bSigma_g(\sS_g\cup \sB_g, \sS_g\cup \sB_g)$ is given by the following four
blocks:
\begin{align*}
  \bSigma_g(\sS_g, \sS_g) & =
  \begin{pmatrix}
    \bSigma_i(\sS_g \cap \sB_i, \sS_g \cap \sB_i) &\bSigma(\sS_g \cap \sB_i, \sS_g \cap \sB_j)
                \\
    \bSigma(\sS_g \cap \sB_j, \sS_g \cap \sB_i) &\bSigma_j(\sS_g \cap \sB_j, \sS_g \cap \sB_j)
  \end{pmatrix}, \\
 \bSigma_g(\sS_g, \sB_g) & =
  \begin{pmatrix}
    \bSigma_i(\sS_g \cap \sB_i, \sB_g \cap \sB_i) &\b 0                               \\
    \b 0 &\bSigma_j(\sS_g \cap \sB_j, \sB_g \cap \sB_j)
  \end{pmatrix},    \\
  \bSigma_g(\sB_g, \sS_g) & =
  \begin{pmatrix}
    \bSigma_i(\sB_g \cap \sB_i, \sS_g \cap \sB_i) &\b 0                               \\
    \b 0 &\bSigma_j(\sB_g \cap \sB_j, \sS_g \cap \sB_j)
  \end{pmatrix},
\end{align*}
and
\begin{align*}
  \bSigma_g(\sB_g, \sB_g) & =
  \begin{pmatrix}
    \bSigma_i(\sB_g \cap \sB_i, \sB_g \cap \sB_i) &\bSigma(\sB_g \cap \sB_i, \sB_g \cap \sB_j)
                \\
    \bSigma(\sB_g \cap \sB_j, \sB_g \cap \sB_i) &\bSigma_j(\sB_g \cap \sB_j, \sB_g \cap \sB_j)
  \end{pmatrix}.
\end{align*}
To compute $\bG^<$, we start from the given $\bA$ and $\bSigma$ and keep applying the update rule \eqref{eq:gless-update-rule}.
Similar to the basic FIND algorithm, if $\sC_g$ is a leaf node, then we have
$\bSigma_g(\sB_g\cup \sS_g, \sB_g \cup \sS_g) = \bSigma(\sB_g\cup \sS_g, \sB_g\cup \sS_g)$, i.e., we can use the entries in the
original $\bSigma$ (see appendix). The recursive process of updating $\bSigma$ is the same as updating $\bA$ in the basic FIND algorithm, but
the update rule is different here.

The pseudocode of the algorithm is an extension of the basic FIND algorithm. In addition to rearranging $\bA$ and computing
matrices $\cU_g$, we need to perform similar operations for $\bSigma$ and compute matrices $\cR_g$. Note that even if we do not
need $\bA^{-1}$, we still need to keep track of the update of $\bA$, i.e., the $\cU$ matrices. This is because we need
$\bA_g(\sB_g\cup \sS_g, \sB_g\cup \sS_g)$ to update $\bSigma_g(\sB_g, \sB_g)$ and obtain $\cR_g$. Once we have prepared $\cR_g$
for all the positive clusters, we can compute $\cR_g$ for the negative clusters. This is done in the downward pass as shown in
the procedure \texttt{updateAdjacentNodes} below. The whole algorithm is shown in procedure \texttt{computeG$^<$}. In these
procedures, we use slightly different notations. Instead of using $\sC_g$, $\sB_g$, and $\sS_g$ for the sets of boundary nodes
and private inner nodes, we use $\sC$, $\sB_{\sC}$, and $\sS_{\sC}$. We also use $\cU_{\sC}$ and $\cR_{\sC}$ for the
corresponding $\cU_g$ and $\cR_g$, respectively.

\LinesNumbered
\begin{procedure}[htbp]
  \caption{updateBoundaryNodes(cluster $\sC$). This procedure is
    called with the root of the tree:
    updateBoundaryNodes(root). \label{algo1}}
  \KwData{tree decomposition of the mesh; the matrix $\bA$.}
  \KwIn{cluster $\sC$ with $n$ boundary mesh nodes.}
  \KwOut{all the inner mesh nodes of cluster $\sC$ are eliminated by the procedure.
    The $n\times n$ matrices $\cU_\sC$ for $\bG^r$ and $\cR_\sC$ for $\bG^<$ are computed and saved.}
  \If{$\sC$ is not a leaf}
  {$\sC$1 = left child of $\sC$\;
    $\sC$2 = right child of $\sC$\;
    updateBoundaryNodes($\sC$1)\tcc*[f]{The boundary set is denoted $\sB_{\sC1}$}\;
    updateBoundaryNodes($\sC$2)\tcc*[f]{The boundary set is denoted $\sB_{\sC2}$}\;
  }
  \Else{
    $\bA_\sC = \bA(\sC, \sC)$\;
  }
  \If{$\sC$ is not the root}
  {
    $\bA_\sC$ = [$\cU_{\sC1}$\ $\bA(\sB_{\sC1}$, $\sB_{\sC2}$);
    $\bA(\sB_{\sC2}$, $\sB_{\sC1}$)\ $\cU_{\sC2}$]\;
    \tcc{A(B$_\text{C1}$,B$_\text{C2}$) and A(B$_\text{C2}$,B$_\text{C1}$) are values from the original matrix A.}
    Rearrange $\bA_\sC$ such that the inner nodes of $\sB_{\sC1} \cup \sB_{\sC2}$ appear first\;
    Set $\sS_\text{g}$ in Eq.~\eqref{eq:gless-update-rule} to be the above inner nodes and $\sB_\text{g}$ the rest\;
    Eliminate the inner nodes to compute $\cU_\sC$ by Eq.~\eqref{eq:major-update}\;
    $\bSigma_\sC$ = [$\bSigma_{\sC1}$ $\bSigma$($\sB_{\sC1}$, $\sB_{\sC2}$);
    $\bSigma$($\sB_{\sC2}$, $\sB_{\sC1}$) $\bSigma_{\sC2}$]\;
    \tcc{$\Sigma$(B$_\text{C1}$, B$_\text{C2}$)
      and $\Sigma$(B$_\text{C2}$, B$_\text{C1}$)
      are values from the original matrix $\Sigma$.}
    Rearrange $\bSigma_\sC$ in the same way as $\bA_\sC$\;
    Compute $\cR_\sC$ based on Eq.~\eqref{eq:gless-update-rule}\;
    Save $\cU_\sC$\ and $\cR_\sC$\;
  }
\end{procedure}

\RestyleAlgo{ruled} \LinesNumbered
\begin{procedure}[htbp]
  \caption{updateAdjacentNodes(cluster $\sC$). This procedure is called with
the root of the tree: updateAdjacentNodes(root).}
  \KwData{tree decomposition of the mesh; the matrix $\bA$; the upward pass [updateBoundaryNodes()] should have been completed.}
  \KwIn{cluster $\sC$ with $n$ adjacent mesh nodes (as the boundary nodes of $\bar{\sC}$).}%

  \KwOut{all the outer mesh nodes of cluster $\sC$ (as the inner nodes of $\bar{\sC}$) are eliminated by the procedure.
  The $n\times n$ matrices $\cU_{\bar{\sC}}$ and $\cR_{\bar{\sC}}$ are computed and saved.}

  \If{$\sC$ is not the root}
  {$\sD$ = parent of $\sC$ \tcc*[f]{The boundary set of $\bar{\text{D}}$ is denoted
      B$_{\bar{\text{D}}}$}\;%
    $\sD$1 = sibling of $\sC$ \tcc*[f]{The boundary set of D1 is denoted
      B$_\text{D1}$}\;%
    $\bA_{\bar{\sC}}$ = [$\cU_{\bar{\sD}}$
    $\bA(\sB_{\bar{\sD}}, \sB_{\sD1})$;
    $\bA(\sB_{\sD1}, \sB_{\bar{\sD}})$ $\cU_{\sD1}$]\;%
    \tcc{A(B$_{\bar{\text{D}}}$,B$_\text{D1}$)
      and A(B$_\text{D1}$,B$_{\bar{\text{D}}}$)
      are values from the original matrix A.}
    \tcc{If D is the root, then ${\bar{\text{D}} = \emptyset}$ and
        A$_{\bar{\text{C}}}$ = R$_{\text{D1}}$.}

    Rearrange $\bA_{\bar{\sC}}$ such that the inner nodes of
$\sB_{\bar{\sD}} \cup \sB_{\sD1}$ appear first\;
    Set $\sS_\text{g}$ in Eq.~\eqref{eq:gless-update-rule} to be the above inner nodes and $\sB_\text{g}$ the rest\;

    Eliminate the outer nodes to compute $\cU_{\bar{\sC}}$ by Eq.~\eqref{eq:major-update}\;

    $\bSigma_{\bar{\sC}}$ = [$\cR_{\bar{\sD}}$
$\bSigma(\sB_{\bar{\sD}}, \sB_{\sD1})$;
    $\bSigma(\sB_{\sD1}, \sB_{\bar{\sD}})$ $\cR_{\sD1}$]\;

    Rearrange $\bSigma_{\bar{\sC}}$ in the same way as for $\bA_{\bar{\sC}}$\;

    Compute $\cR_{\bar{\sC}}$ based on Eq.\eqref{eq:gless-update-rule}\;
    Save $\cU_{\bar{\sC}}$ and $\cR_{\bar{\sC}}$\;
  }

  \If{$\sC$ is not a leaf}
  {
        $\sC$1 = left child of $\sC$\;%
        $\sC$2 = right child of $\sC$\;%
        updateAdjacentNodes($\sC$1)\;%
        updateAdjacentNodes($\sC$2)\;%
  }
\end{procedure}

\LinesNumbered
\begin{procedure}[htbp]
  \caption{computeG$^<$(mesh $\sM$). This procedure is called by any
function that needs $\bG^<$ of the whole mesh.}
  \KwIn{the mesh $\sM$; the matrix $\bA$; the matrix $\bSigma$.}%
  \KwOut{the diagonal entries of $\bG^<$.}%
  Prepare the tree decomposition of the whole mesh\;%
  updateBoundaryNodes(root)\;%
  updateAdjacentNodes(root)\;%
  \For{each leaf node $\sC$}%
  {
        Compute [$\bA^{-1}$]($\sC$, $\sC$) using Eq.~\eqref{eq:cd1}\;%
        Compute $\bG^<_\sC$ based on Eq.~\eqref{eq:gless-nn}\;%
        Save $\bG^<_\sC$ together with its indices\;%
  }%
    Collect $\bG^<_\sC$ with their indices and output all the diagonal entries of $\bG^<$\;
\end{procedure}

\subsubsection{Computation and storage cost}

Similar to the computation cost analysis in \cite{li08}, the major computation cost comes from computing
\eqref{eq:gless-update-rule}.

The cost depends on the size of $\sS_g$ and $\sB_g$. Let $\Ss=|\sS_g|$ and $\Sb = |\sB_g|$ be the sizes,\label{def:size} then the
computation cost for $\cL_g\bSigma_g(\sS_g, \sB_g)$ is $(\frac{1}{3}\Ss^3 + \Ss^2\Sb + \Ss\Sb^2)$ flops. Since $\cL_g$ is already
given in~\eqref{eq:major-update}, the cost is reduced to $\Ss\Sb^2$. Similarly, the cost for $\bSigma_g(\sB_g,
\sS_g)\cL_g^{\dagger}$ is also $\Ss\Sb^2$ and the cost for $\cL_g\bSigma_g(\sS_g, \sS_g)\cL_g^{\dagger}$ is $\Ss\Sb^2 +
\Ss^2\Sb$. So the total cost for \eqref{eq:gless-update-rule} is $(3\Ss\Sb^2 + \Ss^2\Sb)$ flops. Note that these cost estimates
need explicit form of $\bA(\sB, \sS)\bA(\sS, \sS)^{-1}$ and we need to arrange the order of computation to have it as an
intermediate result in computing \eqref{eq:major-update}.

Now we consider the whole process. We analyze the cost in the two passes.

For the upward pass, when two $(a\times a)$-clusters merge into one $(a\times 2a)$-cluster, $\Sb \leq 6a$ and $\Ss \leq 2a$, so
the cost is at most $360a^3$ flops; when two $(a\times 2a)$-clusters merger into one $(2a\times 2a)$-cluster, $\Sb \leq 8a$ and
$\Ss \leq 4a$, so the cost is at most $896a^3$ flops. We have $\frac{N_xN_y}{2a^2}$ $(a\times 2a)$-clusters and
$\frac{N_xN_y}{4a^2}$ $(2a\times 2a)$-clusters, so the cost for these two levels is $(360a^3\frac{N_xN_y}{2a^2} +
896a^3\frac{N_xN_y}{4a^2} =) 404N_xN_ya$. Sum all the levels together up to $a = N_x/2$, the computation cost for upward pass is
then $404N_x^2N_y$.

For the downward pass, when one $(2a\times 2a)$-cluster is partitioned to two $(a\times 2a)$-clusters, $\Sb\leq 6a$ and $\Ss\leq
8a$, so the cost is at most $1248a^3$ flops; when one $(a\times 2a)$-cluster is partitioned to two $(a\times a)$-clusters, $\Sb
\leq 4a$ and $|\sS_g|\leq 6a$, so the cost is at most $432a^3$ flops. We have $\frac{N_xN_y}{2a^2}$ $(a\times 2a)$-clusters and
$\frac{N_xN_y}{a^2}$ $(a\times a)$-clusters, so the cost for these two levels is ($1248a^3\frac{N_xN_y}{2a^2} +
432a^3\frac{N_xN_y}{a^2}=$) $1056N_xN_ya$. Sum all the levels together up to $a = N_x/2$, the computation cost for upward pass is
then at most $1056N_x^2N_y$.

\subsection{Computing off-diagonal entries of $\bG^r$ and $\bG^<$}

\label{sec:current} In addition to computing the diagonal entries of $\bG^r$ and $\bG^<$, the algorithm can also be easily
extended to computing the entries in $\bG^r$ and $\bG^<$ corresponding to neighboring nodes. These entries are often of interest
in simulation. For example, the current density can be computed using off-diagonal entries corresponding to horizontally
neighboring nodes. These entries can be obtained with a slight modification of the algorithm. For simplicity, we only consider
the off-diagonal entries of $\bG^r$. The off-diagonal entries of $\bG^<$ can be treated in a similar way.

\begin{figure}[htbp]
\centering
\subfigcapskip=-28.5pt
\subfigure[t][A leaf cluster with its neighboring nodes. The arrows correspond to the entries for the current density that are being computed.]{%
    \parbox[b]{150pt}{
            \vspace{-30pt}            \hspace*{30pt}\includegraphics[width=95pt]{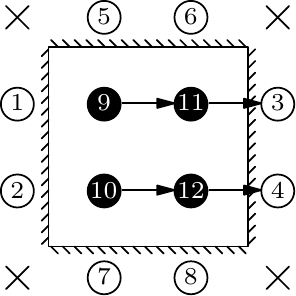}
            \vspace{30pt}
    }
}%
\subfigcapskip=0pt
\hspace{30pt}
\subfigure[Matrix entries.]{%
    \includegraphics[width=150pt]{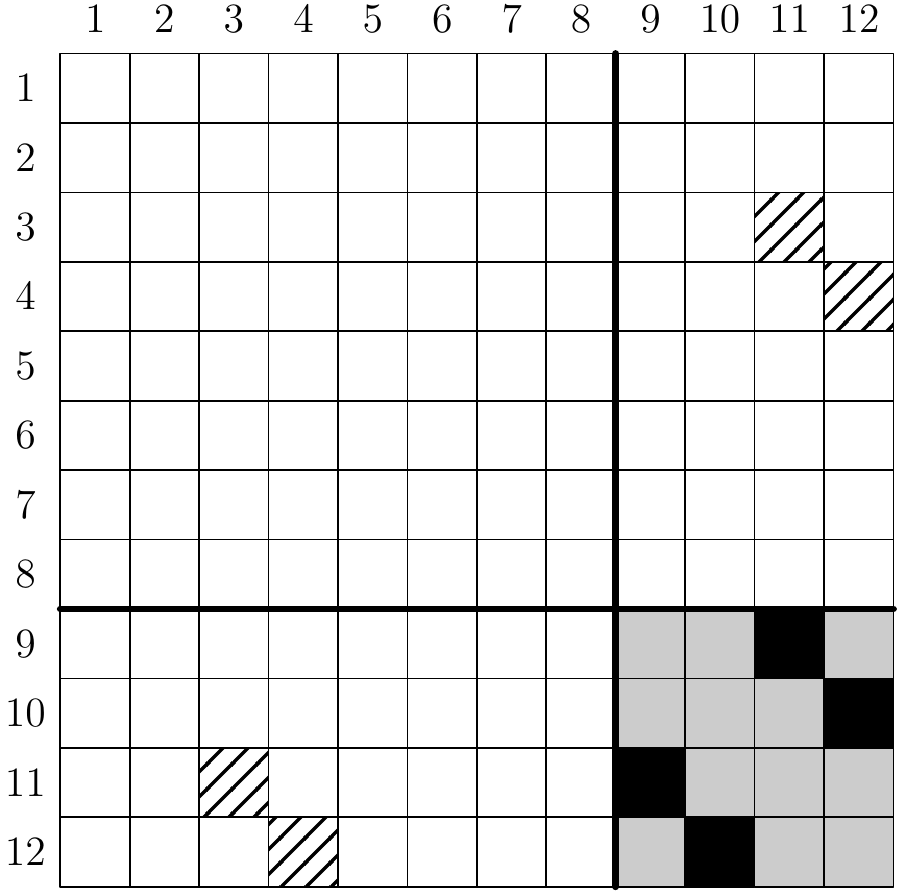}
}
\caption{Last step of elimination. The entries in the shaded area are obtained through computing $\bG^r$.
The solid blocks and the patterned blocks are the off-diagonal entries for the current density.
The solid blocks are obtained for free when we compute the diagonal entries of $\bG^r$. The patterned blocks need to be computed separately.
\label{fig:last-elim}}
\end{figure}

Fig.~\ref{fig:last-elim} shows a small number of mesh nodes and their corresponding matrix entries for the computation of
the diagonal entries $\bG^r$ and the current density. The entries of $\bG^r=\bA^{-1}$ corresponding to the nodes in $\sC$ are
being computed with the following equation:
\begin{equation}
[\bA^{-1}](\sC, \sC) = [\bA(\sC, \sC) - \bA(\sC, \sB_{\bar{\sC}}) (\cU_{\bar{\sC}})^{-1}\bA(\sB_{\bar{\sC}}, \sC)]^{-1}.
\label{eq:cd1}
\end{equation}
The current density between nodes 9 and 11, and that between nodes 10 and 12 are also given by \eqref{eq:cd1}. To compute the
current density between nodes 11 and 3 and that between nodes 12 and 4, perform one step of back substitution and we have
\begin{align}
[\bA^{-1}](\sB_{\bar{\sC}}, \sC) & = -(\cU_{\bar{\sC}})^{-1} \bA(\sB_{\bar{\sC}}, \sC) [\bA^{-1}](\sC, \sC)
\label{eq:cd2} \\
\intertext{and}
[\bA^{-1}](\sC, \sB_{\bar{\sC}}) & = -[\bA^{-1}](\sC, \sC) \bA(\sC, \sB_{\bar{\sC}}) (\cU_{\bar{\sC}})^{-1}.
\label{eq:cd3}
\end{align}
This can be generalized to nodes that are farther away.

\section{Optimization of FIND}\label{sec:optimization}

In our FIND algorithm, we achieved $\cO(N_x^2N_y)$ computation complexity for a 2D mesh of size $N_x\times N_y$. Even though this
complexity has been reduced by an order of magnitude compared to the state-of-the-art RGF method, the matrix inversion is still
the most time consuming part in transport problem simulations. This chapter discusses how to further reduce the computational
cost by using certain properties of the matrix $\bA$ such as its sparsity pattern.

In the FIND algorithm, the major operation is performing Gaussian eliminations. All such operations are of the form
\begin{equation}
  \cU = \bA(\sB, \sB) - \bA(\sB, \sS)\bA(\sS, \sS)^{-1}\bA(\sS, \sB)
  \label{eq:major-update-simple}
\end{equation}
for $\bG^r$ and
\begin{equation}
\cR = \bSigma(\sB, \sB) - \cL\bSigma(\sS, \sB) - \bSigma(\sB, \sS)\cL^{\dagger} + \cL\bSigma(\sS, \sS)\cL^{\dagger}
\label{eq:gless-update-rule-simple}
\end{equation}
for $\bG^<$, where
\begin{equation}
\cL = \bA(\sB, \sS)\bA(\sS, \sS)^{-1}.
\label{eq:express-L}
\end{equation}
These equations are copied from \eqref{eq:major-update} and \eqref{eq:gless-update-rule} with subscripts skipped for
simplicity.\footnote{We use loose notations here since some of the entries of $\bA(\bS, \sS)$ and $\bSigma(\sS, \sS)$ have been
modified by the Gaussian elimination compared to the original $\bA$ and $\bSigma$. See Section~\ref{extension-diagonal} for
details.} As a result, computing \eqref{eq:major-update-simple} and \eqref{eq:gless-update-rule-simple} efficiently becomes
critical for our algorithm to achieve good overall performance.

Simple analysis~\cite{li08} shows that the computation cost for \eqref{eq:major-update-simple} is
\begin{equation}
\frac{1}{3}\Ss^3 + \Ss^2\Sb + \Ss\Sb^2 \label{eq:basic-computation-cost}
\end{equation}
and for \eqref{eq:gless-update-rule-simple} is
\begin{equation}
\Ss^2\Sb + 3\Ss\Sb^2,    \label{eq:gless-computation-cost}
\end{equation}
where $\Ss$ and $\Sb$ are the size of sets $\sS$ and $\sB$, respectively. These costs assume that all the matrices in
\eqref{eq:major-update-simple} and \eqref{eq:gless-update-rule-simple} are general dense matrices. These matrices, however, are
themselves sparse for a typical 2D mesh. In addition, due to the characteristics of the physical problem in many real
applications, the matrices $\bA$ and $\bSigma$ often have special properties~\cite{svizhenko02}. Such sparsities and properties
will not change the order of cost, but may reduce the constant factor of the cost, thus achieving some measure of speed-up. With
proper optimization, the FIND algorithm can exceed other algorithms on a much smaller mesh.

In this section, we first exploit the sparsity of the matrices in \eqref{eq:major-update} and \eqref{eq:gless-update-rule}
to reduce the constant factor in the computation and the storage cost. We then consider the symmetry and positive definiteness of
the given matrix $\bA$ for further performance improvement. Finally, we also apply these optimization techniques to $\bG^<$ and
current density, when $\bSigma$ has similar properties.

\subsection{Optimizations that depend on the sparsity of $\bA$} \label{extra-sparsity}

\label{sec:extra}

Because of the local connectivity of the 2D mesh in our approach\footnote{The approach works for 3D mesh as well}, the matrices
$\bA(\sB, \sS)$, $\bA(\sS, \sS)$, and $\bA(\sS, \sB)$ in \eqref{eq:major-update-simple} are not dense. Such sparsity will not
reduce the storage cost since the matrix $\cU$ in \eqref{eq:major-update-simple} is dense and this is the major part of the
storage cost. However, the computational cost can be significantly reduced. To clearly observe the sparsity and exploit it, it is
convenient to consider these matrices in blocks.

\subsubsection{Schematic description and analysis of the sparsity pattern} \label{sec:block-property}

\begin{figure}[htbp]
\centering
\includegraphics[width=150pt]{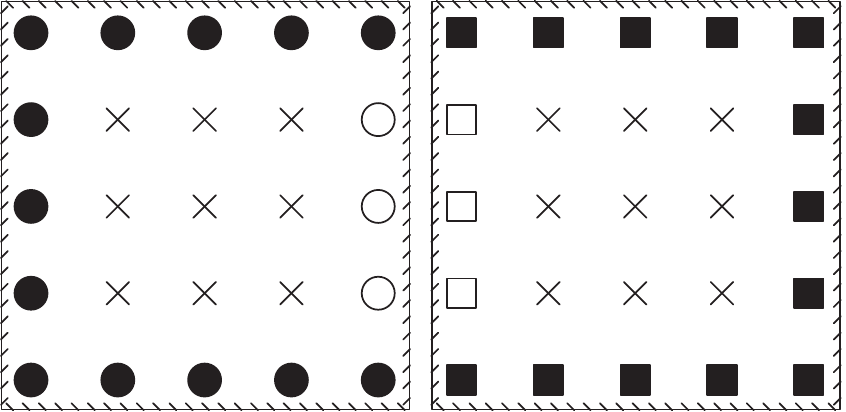}
\caption{ Two clusters merge into one larger cluster in the upward pass. The $\times$ nodes have already been eliminated. The
{\large$\circ$} and {\scriptsize$\square$} nodes remain  to be eliminated. This figure is slightly different from
Fig.~\ref{fig:simple-merge} as we distinguish the nodes in the left and right clusters here.} \label{fig:cluster-merge-up}
\end{figure}

As shown earlier in Fig.~\ref{fig:simple-merge}, we eliminate the private inner nodes when we merge two child clusters in the
tree into their parent cluster. In Figs.~\ref{fig:cluster-merge-up}, we distinguish between the left and right clusters in the
merge (compare with Fig.~\ref{fig:simple-merge}).

In Fig.~\ref{fig:cluster-merge-up}, the three hollow circle nodes and the three hollow square nodes are private inner nodes of
the parent cluster. They originate from the left child cluster and the right child cluster and are denoted as $\sS_L$ and
$\sS_R$,\label{def:SL-SR} respectively. When we merge the two child clusters, these nodes will be eliminated. Similarly, the
solid circle nodes and the solid square nodes are boundary nodes originated from the left and the right child clusters and are
denoted as $\sB_L$ and $\sB_R$,\label{def:BL-BR} respectively. When we merge the two child clusters, these nodes will be updated.
The matrix blocks corresponding to the circle and square nodes can be written as
\begin{equation}
\begin{pmatrix}
\bU(\Circ, \Circ) &\bA(\Circ, \Square)  &\bU(\Circ, \Bullet)    &\b 0  \\
\bA(\Square, \Circ) &\bU(\Square, \Square)  &\b 0     &\bU(\Square, \Blacksquare)   \\
\bU(\Bullet, \Circ)    &\b 0  &\bU(\Bullet, \Bullet) &\bA(\Bullet, \Blacksquare) \\
\b 0  &\bU(\Blacksquare, \Square)    &\bA(\Blacksquare, \Bullet)    &\bU(\Blacksquare, \Blacksquare)
\end{pmatrix},
\label{eq:block-elim-matrix}
\end{equation}
and the merge corresponds to the elimination of the first two columns of \eqref{eq:block-elim-matrix}. If we follow the notation
used in Section~\ref{sec:basic-detailed}, then \eqref{eq:block-elim-matrix} corresponds to $\bA_k(\sS_k\cup\sB_k,
\sS_k\cup\sB_k)$\footnote{We use $k$ instead of $g$ here.}, where cluster $k$ has two child clusters $i$ and $j$, and the circle
and square nodes are given by:
\begin{lsitemizei}
\item[\qquad$\Circ$\,:] $\sS_L = \sS_k\cap \sB_i$, private inner nodes from the left cluster $i$
\item[\qquad$\Square$\,:] $\sS_R = \sS_k\cap \sB_j$, private inner nodes from the right cluster $j$
\item[\qquad$\Bullet$\,:] $\sB_L = \sB_k\cap \sB_i$, boundary nodes from the left cluster $i$
\item[\qquad$\Blacksquare$\,:] $\sB_R = \sB_k\cap \sB_j$, boundary nodes from the right cluster $j$
\end{lsitemizei}

In \eqref{eq:block-elim-matrix} we use $\bU$ instead of $\bA$ for some blocks to emphasize that they are results from earlier
eliminations while the other blocks come from the original matrix $\bA$.\footnote{We use loose notation here. The notations $\bU$
for some of the blocks (e.g., $\bU(\bullet, \bullet)$) is incorrect since this is not the final $\bU$ matrix from the LU
factorization, but rather an intermediate step in the method. This notation is closer to a computer code implementation. More
precise notations however would reduce the legibility of the text.} These $\bU$ blocks are typically dense, since after the
elimination of the inner nodes, the remaining nodes become fully connected. In contrast, the $\bA$ blocks are typically sparse.
For example, since each $\Circ$ node is connected with only one $\Square$ node in Fig.~\ref{fig:cluster-merge-up}, $\bA(\Circ,
\Square)$ and $\bA(\Square, \Circ)$ are both diagonal. $\bA(\Bullet, \Blacksquare)$ and $\bA(\Blacksquare, \Bullet)$ are almost
all 0 except for a few entries, although such sparsity saves little cost because they are only involved in addition operations.

\begin{figure}[htbp]
\centering
\includegraphics[width=240pt]{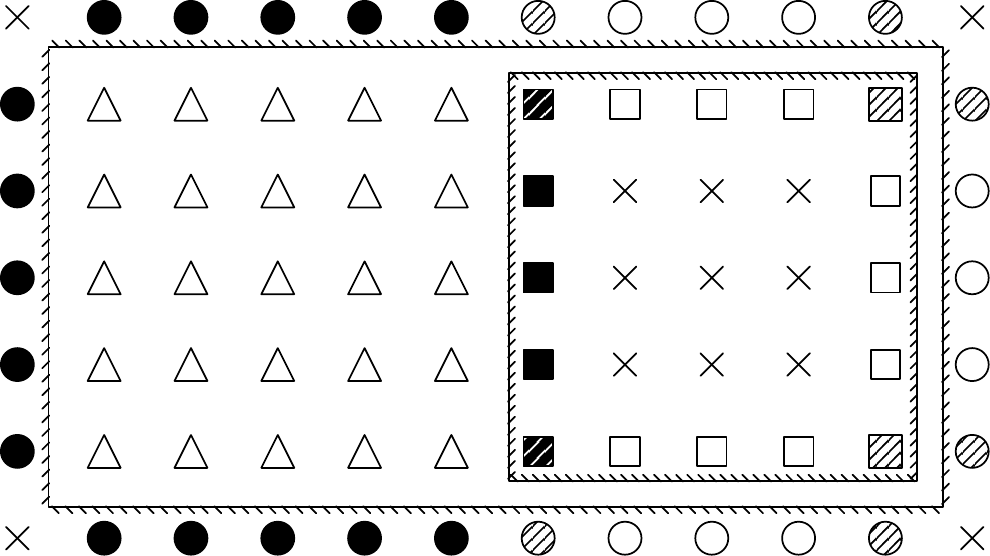}
\caption{Notations are the same as in Fig.~\ref{fig:cluster-merge-up}, with some nodes patterned. The $\vartriangle$ nodes are
not modified in this merging step. The patterned symbols correspond to nodes whose connectivity is different from
Fig.~\ref{fig:cluster-merge-up}. As a result, $\bA$({\large$\circ$}, {\scriptsize$\square$}$)$ and $\bA(${\scriptsize$\square$},
{\large$\circ$}) are no longer strictly diagonal and $\bA(${\large$\circ$}, {\scriptsize$\blacksquare$}$)$ is not exactly zero
any more.}\label{fig:cluster-merge-down}
\end{figure}

These properties always hold for the upward pass. They hold for the downward pass as well but with a few exceptions. These
exceptions are illustrated by the patterned nodes in Fig.~\ref{fig:cluster-merge-down}. Since the exceptions only occur at the
corner nodes, we can ignore them in the computational and storage cost analysis.

\subsubsection{Exploiting the sparsity using a block structure} \label{sec:exploit}

We analyzed the computational cost using different formulas to calculate the inverse of a general $2\times 2$ block matrix. The following forms of the
inverse are available:\footnote{We assume $A$ and $D$ to be nonsingular.}
\begin{subequations}  \label{eq:block-inverse}
\begin{eqnarray}
\begin{pmatrix}
A   &B  \\
C   &D
\end{pmatrix} ^{-1} &= &\begin{pmatrix}
                        \tilde{A}^{-1}              &-A^{-1} B \tilde{D}^{-1}   \label{def:A-D-tilde}   \\
                        -D^{-1} C \tilde{A}^{-1}    &\tilde{D}^{-1}
                        \end{pmatrix}   \label{eq:back-sub-2way}    \\
                    &= &\begin{pmatrix}
                        A^{-1} + A^{-1}B\tilde{D}^{-1}CA^{-1}              &-A^{-1} B \tilde{D}^{-1}    \\
                        -\tilde{D}^{-1} C A^{-1}    &\tilde{D}^{-1}
                        \end{pmatrix}   \label{eq:back-sub-1way}    \\
                    &= &\begin{pmatrix}
                        A   &B  \\
                        0   &\tilde{D}
                        \end{pmatrix}^{-1}
                        \begin{pmatrix}
                        I   &0  \\
                        CA^{-1}   &I
                        \end{pmatrix}^{-1},\label{eq:block-LU-inverse}%
\end{eqnarray}%
\label{eq:inverse-methods}%
\end{subequations}%
where $\tilde{A} = A - BD^{-1}C$ and $\tilde{D} = D - CA^{-1}B$ are the Schur complements of $D$ and $A$, respectively.
In~\eqref{eq:back-sub-2way}, the calculations for $\tilde{A}^{-1}$ and $\tilde{D}^{-1}$ are independent of each other and can be
done in parallel, so we call this method \emph{parallel inverse}. In~\eqref{eq:back-sub-1way}, we have to calculate
$\tilde{D}^{-1}$ first to calculate the other block of the inverse: $A^{-1} + A^{-1}B\tilde{D}^{-1}CA^{-1}$, so we call this
method \emph{sequential inverse}. In~\eqref{eq:block-LU-inverse}, we perform block LU factorization first and then calculate the
inverse, so we call this method \emph{block LU inverse}.

Details of our derivation can be found in the appendix. Using our analysis of the sparsity pattern, the cost for each approach is shown in the table below:
\begin{table}[htbp]
\centering \caption{Summary of the optimization methods.}
\label{operation-cost-summary}
\begin{tabular}{@{}cccc@{}}    \toprule
     \emph{Method}     &\emph{Cost}     &\emph{Reduction in four cases}\footnotemark  &\emph{Parallelism}   \\ \midrule
	\addlinespace[2pt]
    parallel inverse    &$\frac{8}{3}m^3 + 4m^2n + 4mn^2$   &51.4, 52.6, 61.5, 60.1             &good   \\[2pt]
    sequential inverse  &$\frac{4}{3}m^3 + 5m^2n + 4mn^2$      &53.0, 54.0, 55.8, 55.7         &little \\[2pt]
    block LU inverse    &$\frac{4}{3}m^3 + 4m^2n + 5mn^2$       &59.1, 57.9, 53.9, 54.3 &some      \\[2pt]
    na\"{i}ve LU        &$\frac{8}{3}m^3 + \frac{13}{2}m^2n + 4mn^2$    &59.0, 62.5, 76.0, 74.4        & little   \\  \bottomrule
\end{tabular}
\end{table}
\footnotetext{The percentage of the cost given by \eqref{eq:basic-computation-cost}.}

Note that these methods are not exhaustive. There are many variations of these methods with different forms of $\bA(\sS, \sS)^{-1}$, different computation orders, and different common parts in the computation. Determining which method is the best may depend on the size of $\sS_L$, $\sS_R$, $\sB_L$, and $\sB_R$, the exceptions in the downward pass (corner nodes), the flop rate, the cache miss rate, and the implementation complexity.

\subsection{Optimization for symmetry and positive definiteness properties} \label{symmetry}

In real problems, $\bA$ is often symmetric (or Hermitian if they are complex, which can be treated in a similar
way)~\cite{svizhenko02}. So it is worthwhile to consider special treatment for such matrices to reduce both computation cost and
storage cost. Note that this reduction is independent of the optimizations based on the sparsity pattern in Section.~\ref{extra-sparsity}.

If all the given matrices are symmetric, it is reasonable to expect the elimination results to be symmetric as well since our
update rule \eqref{eq:major-update-simple} does not break matrix symmetry. This is shown in the following \emph{property of
symmetry preservation}:
\begin{property}[symmetry preservation]
If $\bA$ is symmetric, then in~\eqref{eq:major-update-simple}, $\cU, \bA(\sB, \sB)$, and $\bA(\sS, \sS)$ are all symmetric;
$\bA(\sB, \sS)$ and $\bA(\sS, \sB)$ are transpose of each other. \label{symmetry-preservation}
\end{property}
See \ref{proof_sympres} for a proof.

In addition, $\bA$ is often positive definite. This property is also preserved during the elimination process in our algorithm,
as shown in the following \emph{property of positive-definiteness preservation}:

\begin{property}[positive-definiteness preservation]
If $\bA$ is symmetric and positive-definite, then the matrix $\cU$ in \eqref{eq:major-update-simple} is always positive-definite.
\label{pd-preservation}
\end{property}
See \ref{proof_posdef} for a proof. Using the symmetry preservation property and the positive definiteness preservation property, we can write the last term in
\eqref{eq:major-update-simple} as
$$\bA(\sB, \sS)\bA(\sS, \sS)^{-1}\bA(\sS, \sB) = (\cG_{\sS}^{-1}\bA(\sS, \sB))^T (\cG_{\sS}^{-1}\bA(\sS, \sB)),$$
where $\bA(\sS, \sS) = \cG_{\sS}\cG_{\sS}^T$ is the Cholesky factorization of $\bA(\sS, \sS)$. The Cholesky factorization has
cost $\frac{1}{6}\Ss^3$. The solver has cost $\frac{1}{2}\Ss^2\Sb$, and the final multiplication has cost $\frac{1}{2}\Ss\Sb^2$.
The cost for \eqref{eq:major-update-simple} is then reduced by half from $\frac{1}{3}\Ss^3 + \Ss^2\Sb + \Ss\Sb^2$ to
$\frac{1}{6}\Ss^3 + \frac{1}{2}\Ss^2\Sb + \frac{1}{2}\Ss\Sb^2$.

Note that even if $\bA$ is not positive definite, by the symmetry preservation property alone, we can still write the last term
of \eqref{eq:major-update-simple} as
$$\bA(\sB, \sS) \bA(\sS, \sS)^{-1} \bA(\sS, \sB) = (\cL_{\sS}^{-1}\bA(\sS, \sB))^T \bD_{\sS}^{-1} (\cL_{\sS}^{-1}\bA(\sS, \sB)),$$
where $\bA(\sS, \sS) = \cL_{\sS}\bD_{\sS}\cL_{\sS}^T$ is the LDL$^{\mbox{\scriptsize T}}$ factorization of $\bA(\sS, \sS)$.

Similarly, the computation cost is reduced by half. However, this method may be subject to large errors due to small pivots,
since the pivots can only be selected from the diagonal entries~\cite{golub96}. The diagonal pivoting
method~\cite{bunch71,bunch77} can be used to solve such a problem, but it is computationally more expensive.

\paragraph{Combining symmetry and sparsity}
For symmetric and positive definite $\bA$, we can take advantage of the sparsity pattern discussed in Section~\ref{extra-sparsity}
as well. Consider the following block Cholesky factorization of $\bA(\sS, \sS)$:
\begin{equation} \label{eq:symmetry-A1}
\bA(\sS, \sS) =
\begin{pmatrix}
    \bA_{11}  &\bA_{12} \\
    \bA_{12}^T  &\bA_{22}
\end{pmatrix}
=
\begin{pmatrix}
    \cG_1   &0  \\
    \bA_{12}^T\cG_1^{-T} &\widetilde{\cG_2}
\end{pmatrix}
\begin{pmatrix}
    \cG_1^T &\cG_1^{-1}\bA_{12} \\
    0   &\widetilde{\cG_2}^T
\end{pmatrix} = \cG_{\sS}\cG_{\sS}^T,
\end{equation}
where $\bA_{11} = \cG_1\cG_1^T, \widetilde{\bA_{22}} = \bA_{22} - \bA_{12}^T\bA_{11}^{-1}\bA_{12}$, and $\widetilde{\bA_{22}} =
\widetilde{\cG_2}\widetilde{\cG_2}^T$. Now, we have
\begin{eqnarray}    \label{eq:symmetry-A2}
\cG_{\sS}^{-1} \bA(\sS, \sB)
    &= &\blockm{\cG_1}{0}{\bA_{12}^T\cG_1^{-T}}{\widetilde{\cG_2}}
        \blockm{\bA_i(\sS_k\cap \sB_i, \sB_k\cap \sB_i)}{0}{0}{\bA_j(\sS_k\cap \sB_j, \sB_k\cap \sB_j)} \notag \\
    &= &\blockm{\cG_1^{-1}\bA(\sS_L, \sB_L)}{0}
        {-\widetilde{\cG_2}^{-1}\bA_{12}^T \cG_1^{-T} \cG_1^{-1}\bA(\sS_L, \sB_L)}{\widetilde{\cG_2}^{-1} \bA(\sS_R, \sB_R)}.
\end{eqnarray}
and
\begin{equation}    \label{eq:symmetry-A3}
\bA(\sB, \sS) \bA(\sS, \sS)^{-1} \bA(\sS, \sB) = (\cG_{\sS}^{-1}\bA(\sS, \sB))^T (\cG_{\sS}^{-1}\bA(\sS, \sB)).
\end{equation}
With these different optimizations, the total cost ends up being $\frac{2}{3}m^3 + 2m^2n + \frac{5}{2}mn^2$.
Compared to the original cost with a block structure, the cost is reduced by more than half when $m=n$. The cost with both optimizations, for sparsity and symmetry, is 25\% of the original cost when the condition $m \approx \frac{\sqrt{3}}{2}n$ is satisfied.


\paragraph{Storage cost reduction}

In addition to the computation cost reduction, the storage cost can also be reduced for symmetric $\bA$. Since the storage cost
is mainly for $\cU$, which is always symmetric, the storage cost is reduced by half. Another part of the storage cost comes from
the temporary space needed for \eqref{eq:major-update-simple}, which can also be reduced by half. Such space is needed for the
top level update in the cluster tree. It is about the same as the cost for $\cU$ but we do not have to keep it. So for a typical
cluster tree with ten or more levels, this part is not important. When the mesh size is small with a short cluster tree, this
part and its reduction may be important.


\subsection{Optimization for computing $\bG^<$ and current density}

We can also reduce the cost for computing $\bG^<$ by optimizing the procedure for \eqref{eq:gless-update-rule-simple}:
\begin{equation}
\cR = \bSigma(\sB, \sB) - \cL\bSigma(\sS, \sB) - \bSigma(\sB, \sS)\cL^{\dagger} + \cL\bSigma(\sS, \sS)\cL^{\dagger}
\label{eq:gless-update-rule-simple2}
\end{equation}
We now consider optimizations that use the sparsity pattern of $\bA$ and $\bSigma$, and separately the symmetry of $\bA$ and $\bSigma$. Those optimizations can be combined but this will not be presented in this paper.

\subsubsection{$\bG^<$ sparsity}

Since $\bSigma(\sS, \sB)$ and $\bSigma(\sB, \sS)$ are block diagonal, the cost for the second and the third term in
\eqref{eq:gless-update-rule-simple2} is reduced by half. Similar to the structure of $\bA(\sS, \sS)$ in
\eqref{eq:major-update-simple}, the blocks $\bSigma(\sS_k \cap \sB_i, \sS_k \cap \sB_j)$ and $\bSigma(\sS_k \cap \sB_j, \sS_k
\cap \sB_i)$ in $\bSigma(\sS_k, \sS_k)$ are diagonal. Consequently, the computational cost of operating with these blocks will be
ignored since it is negligible. The cost for the fourth term $\cL^{-1}\bSigma(\sS, \sS)\cL^{-\dagger}$ in
\eqref{eq:gless-update-rule-simple2} is then reduced to $\frac{1}{2}\Ss^2\Sb + \Ss\Sb^2$ (or $\Ss^2\Sb + \frac{1}{2}\Ss\Sb^2$,
depending on the order of computation). So the total cost for \eqref{eq:gless-update-rule-simple2} becomes $\frac{1}{2}\Ss^2\Sb +
2\Ss\Sb^2$ (or $\Ss\Sb + \frac{3}{2}\Ss\Sb^2$). Compared with the cost without optimization for sparsity, it is reduced by 37.5\%
when $\Ss=\Sb$.

In addition, $\bSigma$ is often diagonal or 3-diagonal in real problems. As a result, $\bSigma(\sS_k\cap \sB_i, \sS_k\cap \sB_j)$
and $\bSigma(\sS_k\cap \sB_j, \sS_k\cap \sB_i)$ become zero and $\bSigma_k(\sS_k, \sS_k)$ becomes block diagonal. This can lead
to further reduction.


\subsubsection{$\bG^<$ symmetry}\label{sec:gless-symmetry}

For symmetric $\bSigma$, we have a property of symmetry and positive definiteness preservation for $\bSigma_{g+}(\sB_g, \sB_g)$, $\bSigma_g(\sB_g, \sB_g)$, and $\bSigma_g(\sS_g, \sS_g)$. With such symmetry, we can perform a Cholesky factorization on $\bSigma(\sS, \sS) = \cK_{\sS}\cK_{\sS}^T$. The total cost for \eqref{eq:gless-update-rule} is reduced to $\frac{1}{6}\Ss^3 + \Ss^2\Sb + \frac{3}{2}\Ss\Sb^2$. Compared to the cost of computation without exploiting the symmetry ($\Ss^2\Sb + 3\Ss\Sb^2$), the cost is reduced by approximately 33\% when $\Ss=\Sb$.



\subsubsection{Current density}

\begin{figure}[htbp] \centering
	\parbox[c]{90pt}{\includegraphics[width=90pt]{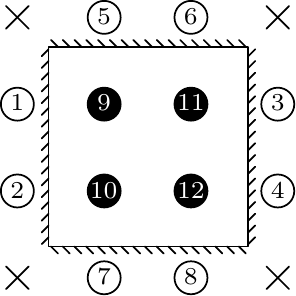}} \hspace{40pt}
    \parbox[c]{140pt}{\includegraphics[width=140pt]{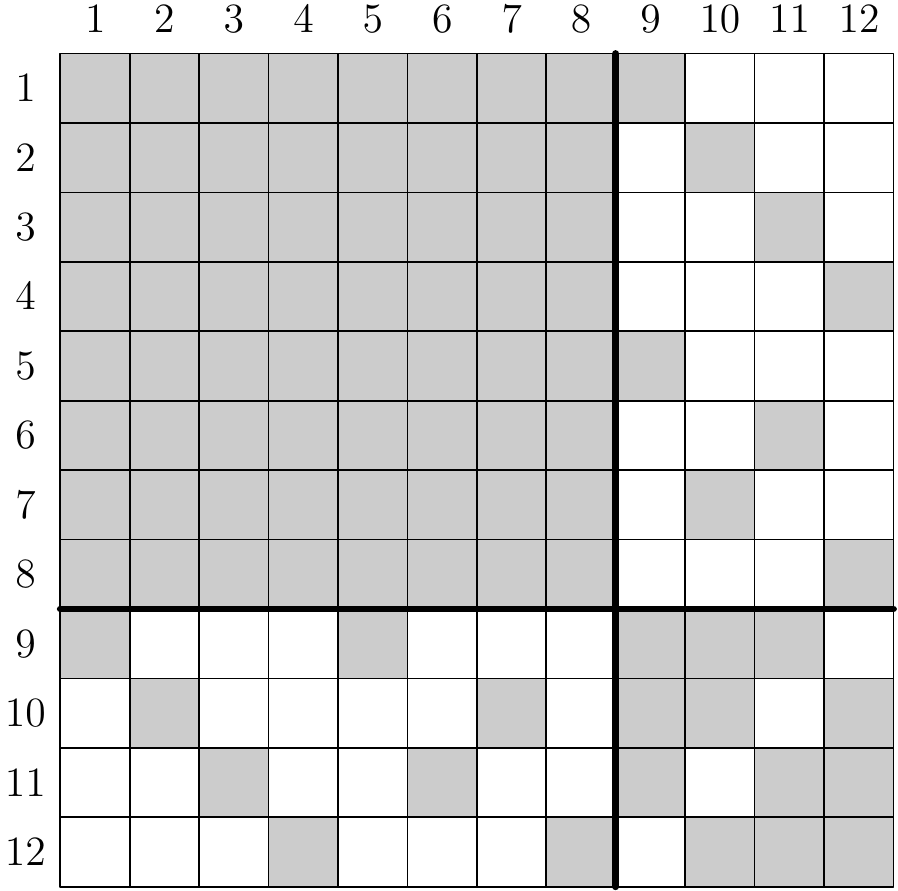}}
\caption{Sparsity of the matrices before the last step of elimination.
Shaded blocks are non-zero entries and blank blocks are zero entries.}
\label{fig:matrix-last-elim-sparsity}
\end{figure}

We can also exploit the sparsity when computing the current density, but it cannot improve much since there is not as much
sparsity in \eqref{eq:cd1}--\eqref{eq:cd3} as in \eqref{eq:major-update-simple}. This can been seen in
Fig.~\ref{fig:matrix-last-elim-sparsity}.

\begin{figure}[htbp]
\centering
\subfigure[Coverage of all the needed current densities using half of the leaf clusters. ]{%
    \centering
    \includegraphics[width=150pt]{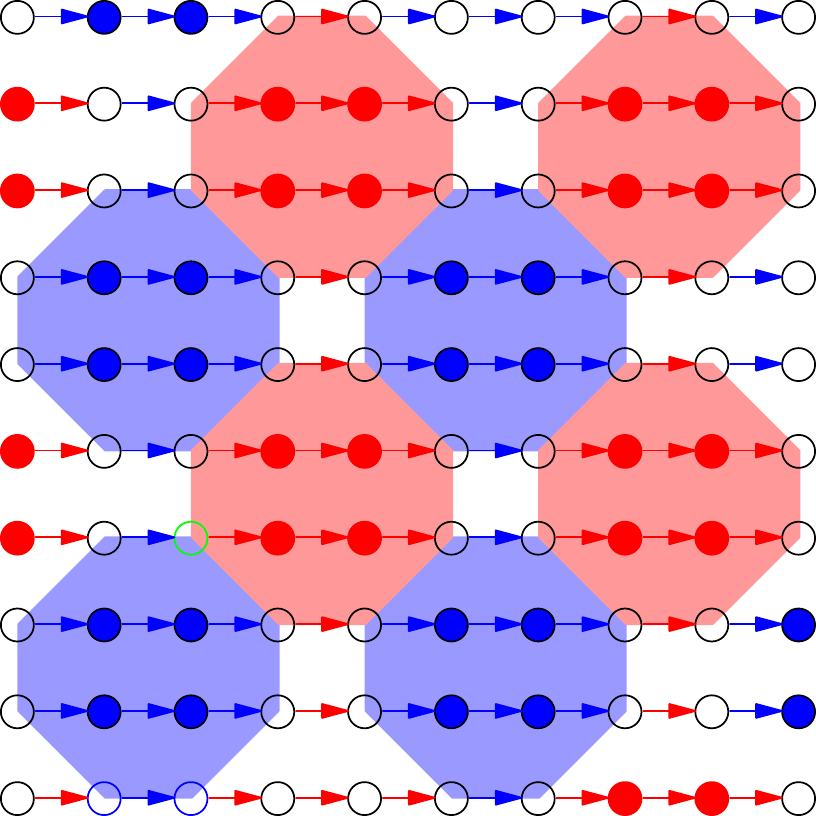} \label{fig:nodes-layout}
}%
\hspace{25pt}
\subfigure[Matrix entries. The notations here follow those in Fig.~\ref{fig:last-elim}.
]{%
    \includegraphics[width=150pt]{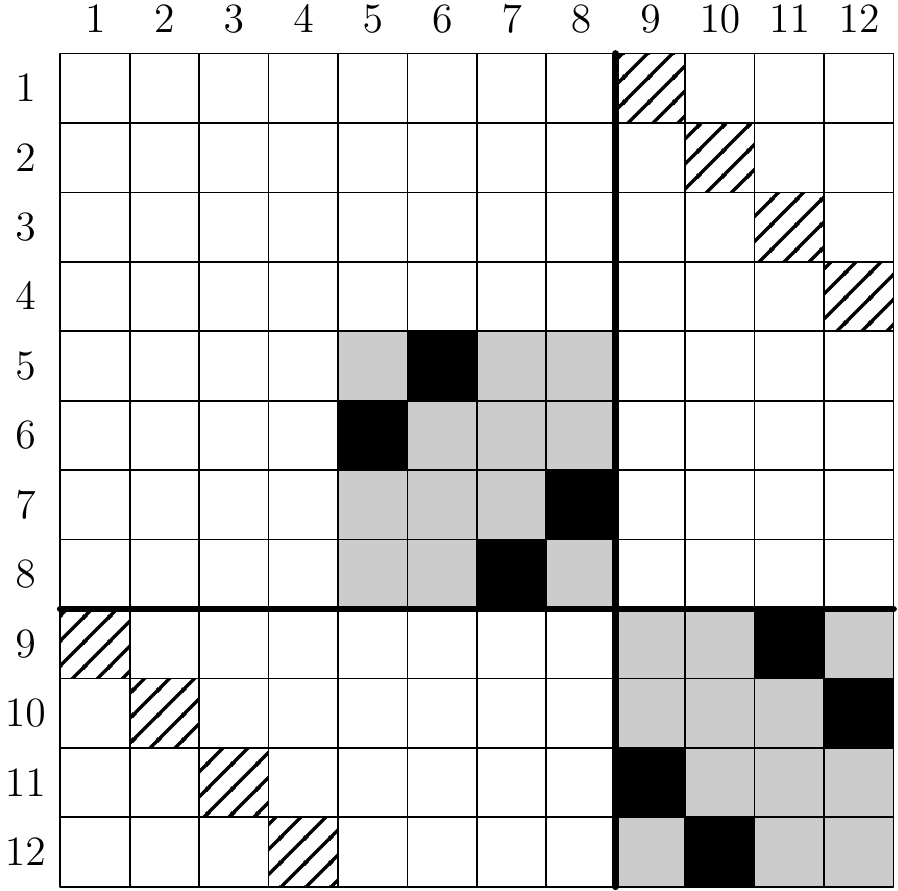}    \label{fig:matrix-entry-details}
}%
\caption{Choice of leaf clusters and the corresponding matrix with desired entries marked.
Each needed cluster consists of four solid nodes. The circle nodes are ``skipped.'' Arrows correspond to the current densities.
If we compute the entries in the shaded area and in the patterned blocks,
then all the desired entries will be covered. See Fig.~\ref{fig:matrix-last-elim-sparsity} for the node numbering.
} \label{fig:cd-opt} \vspace{-20pt}
\end{figure}

However, we can significantly reduce the cost with another approach. In this approach, we still compute the entries of $\bG^r$
for every leaf cluster, but compute the current density only for half of the leaf clusters indicated by the solid nodes in
Fig.~\ref{fig:nodes-layout}. This approach is still based on~\eqref{eq:cd1}--\eqref{eq:cd3}, but in addition to computing the
current density for nodes 9--12, we also compute the current density for nodes 1 and 2, which incurs no extra cost. Therefore, it
reduces the cost for the current density by half.

Moreover, if for the chosen leaf clusters, we compute not only the off-diagonal entries for the current density, but also the
diagonal entries for nodes 5--12, then we can completely avoid the last step of
elimination and the computation thereafter for the other half of the leaf clusters in the final stage. Detailed analysis shows
that computing both the diagonal entries of $\bG^r$ and the current density with this approach is in fact even less costly than computing
only the diagonal entries of $\bG^r$ without this approach. Fig.~\ref{fig:cd-opt} shows how we choose the leaf clusters to cover
the desired diagonal and off-diagonal entries in $\bG^r$. Note that this approach (even though more efficient) was not implemented in our numerical benchmarks.

We can also optimize \eqref{eq:cd1}--\eqref{eq:cd3} by similar techniques used in Section~\ref{symmetry} to reduce the cost for
symmetric $\bA$ and $\bSigma$. In \eqref{eq:cd1},
\begin{equation*}
	\bA(\sC, \sC) - \bA(\sC, \sB_{\bar{\sC}})
	(\cU_{\bar{\sC}})^{-1}\bA(\sB_{\bar{\sC}}, \sC)
\end{equation*}
is of the same form as \eqref{eq:major-update-simple} so the cost is reduced by
half. The cost of computing its inverse is also reduced by half when it is symmetric. Eqs.~\eqref{eq:cd2} and \eqref{eq:cd3} are
transpose of each other so we only need to compute one of them. As a result, the cost for computing the current density is
reduced by half using symmetry.

\subsection{Nodes on the physical boundary of the device} \label{sec:null-boundary}

Some of the nodes on the physical boundary of the device have a different connectivity. For example, a typical node in a 5 point stencil is connected to 4 neighbors. A node on the physical boundary of the device may be connected to only 3 neighbors. This reduces the computational cost of our method since those nodes can be eliminated early in the scheme leading to overall smaller boundary sets for some of the clusters. This has a small effect for leaf clusters (i.e., small clusters) but leads to a significant reduction in the computational cost near the root of the tree when we have few very large clusters. With such effect taken into consideration, the total cost of the algorithm is approximately $457N^3$ (see the appendix). As a comparison, the cost using RGF with symmetry taken into consideration is $3.5N^4$. Since the cross-point of the two cost curves is at $N$ equal to the ratio between the constant factors in the two methods, it is expected to happen at $N\thicksim130$. This is consistent with what we observed in our simulations.

\medskip

\noindent {\it Section Summary.}\ \ The optimization for sparsity leads to $40\%$ cost reduction and that for symmetry leads to $50\%$ cost reduction. This results
in a $70\%$ cost reduction. In other words, the new cost is about \textbf{1/3} of the original cost. This is reflected in the
cross-point of the performance curves using these two methods. Originally, the cross-point was around $N\thicksim130$, but now it can be
as low as $N=40$. This shows the performance improvement resulting from our optimizations. The optimization for sparsity does not
reduce the memory cost. The optimization for symmetry does reduce it (approximately by half).

\section{Numerical Results} \label{sec:num}

\subsection{Extension}

\begin{figure}[htbp] \centering
\subfigure[Running time comparison] {%
    \centering
    \includegraphics[width=150pt]{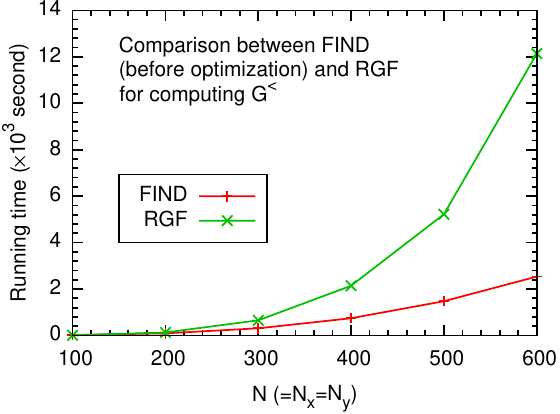}
    \label{fig:Gless-comparison-normal}
}%
\hspace{25pt}%
\subfigure[Running time in log-log scale to see the asymptotic behavior] {%
    \includegraphics[width=150pt]{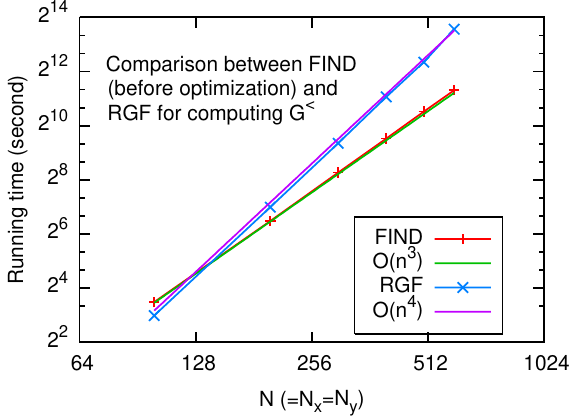} \hspace{5pt}
    \label{fig:Gless-comparison-loglog}
}%
\caption{ Comparison of running time based on FIND extension (before optimization) and RGF for computing $\bG^<$} \label{fig:Gless-comparison}
\end{figure}

In the following two figures, Fig.~\ref{fig:Gless-comparison-normal} gives the running time for computing $\bG^<$ based on the
extension of FIND as described in Section~\ref{sec:extension} and the comparison with RGF. Fig.~\ref{fig:Gless-comparison-loglog}
shows the same data in log-log scale to see more clearly the asymptotic behavior of these two methods. All the running times in
these two figures are based on the extension of FIND before the optimizations discussed in Section~\ref{sec:optimization}.

\subsection{Optimization}

\begin{figure}[htbp]
\centering
\subfigure[Optimization based on sparsity] {%
    \centering
    \includegraphics[width=152pt]{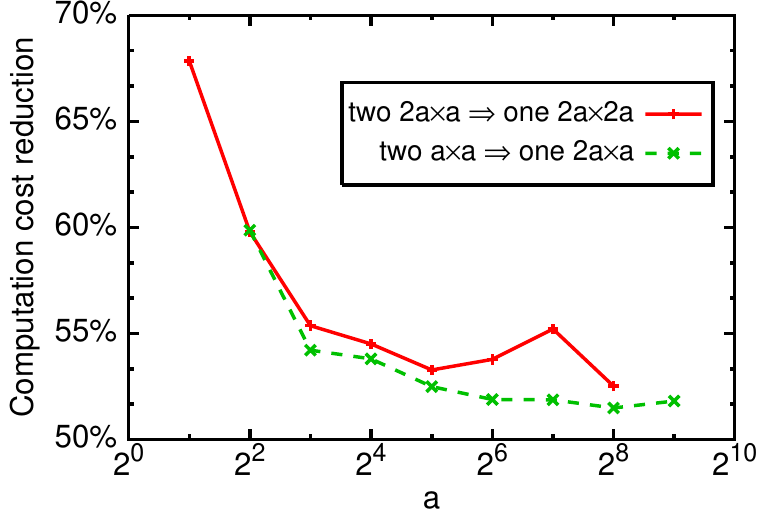}
    \label{extra-sparsity-upward}
}%
\hspace{14pt}%
\subfigure[Optimization based on symmetry] {%
    \includegraphics[width=152pt]{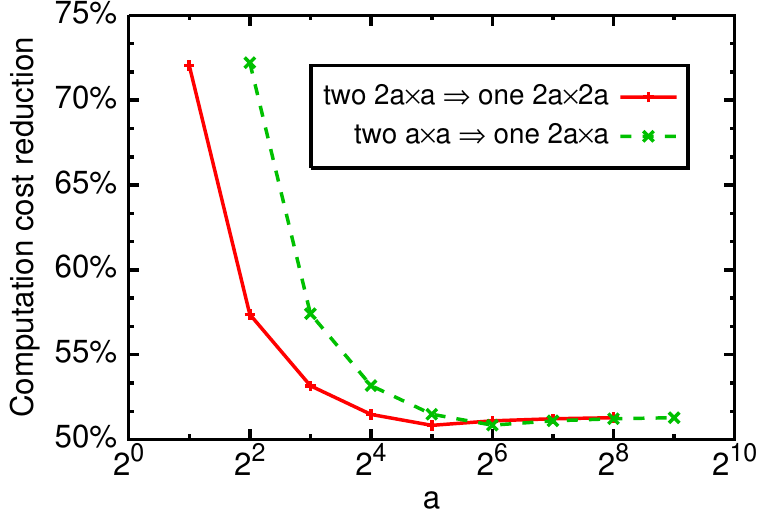}
    \label{optimization-for-symmetry}
}%
\caption{ Comparison between RGF and FIND after optimization\label{fig:comparison-after-optimization}} 
\end{figure}

Fig.~\ref{extra-sparsity-upward} shows the computation cost reduction for clusters of different sizes during the upward pass,
after optimization for sparsity. As indicated in Section~\ref{extra-sparsity}, the improvements are different for two types of
merge: i) two $(a\times a)$-clusters $\Rightarrow$ one $(a\times 2a)$-cluster; ii) two $(2a\times a)$-clusters $\Rightarrow$ one
$(2a\times 2a)$-cluster, so they are shown separately in the figure. The reduction for small clusters is not as significant. This
is because the second order contributions to the cost (e.g., $\cO(m^2)$) was ignored in our analysis but is significant for small
clusters.

Fig.\ref{optimization-for-symmetry} shows the computation cost reduction after optimization for symmetry. Similar to the
Fig.~\ref{extra-sparsity-upward}, we show the reduction when merging clusters at each level in the basic cluster tree. We can see
that the cost is reduced almost by half for clusters larger than $32\times 32$.
\begin{figure}[htbp] \centering
\includegraphics[width=235pt]{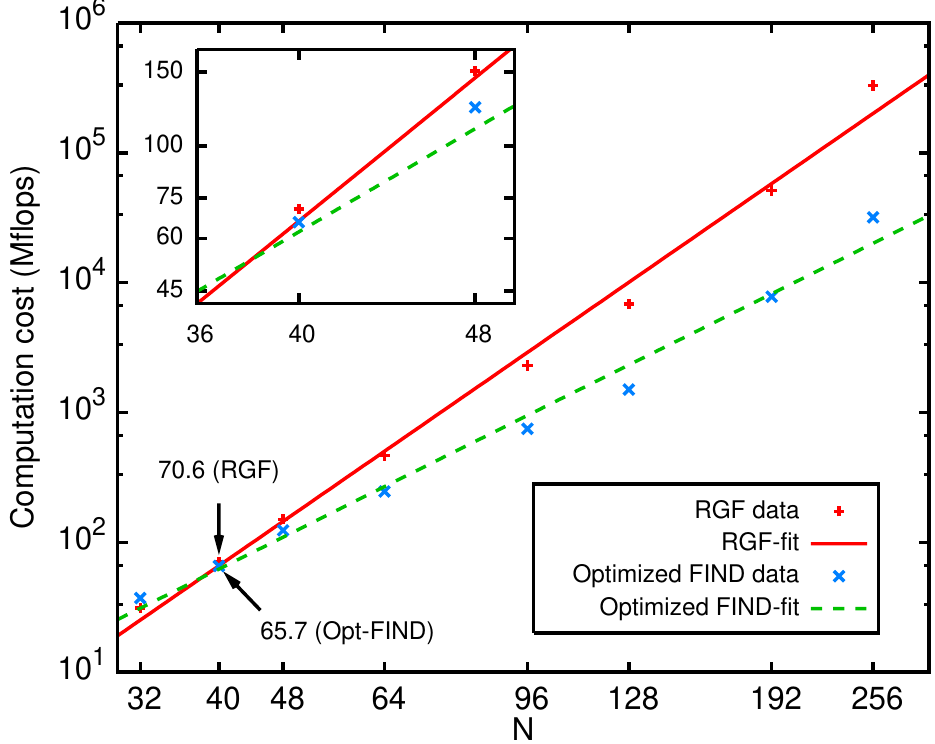} \hspace{5pt}
\caption{ Comparison with RGF to show the change of cross-point} \label{fig:comparison-with-RGF}
\end{figure}
The smaller reduction for small clusters is mainly due to the fact that the Cholesky factorization does not provide savings by a
factor of 2 for these clusters.
Although the reduction in cost is smaller for small clusters using either optimizations, these clusters do not play a significant
role for a mesh larger than $32\times 32$ since the top levels of the tree in Fig.~\ref{fig:cluster-tree-upward} dominate the
computation cost.

Fig.~\ref{fig:comparison-with-RGF} compares the computation cost of RGF and optimized FIND on a Hermitian and positive-definite
matrix $\bA$ and shows the overall effect of the optimization for sparsity and symmetry. The simulations were performed on a
square mesh of size $N\times N$. The cross-point is around $40\times 40$ according to the two fitted lines. We confirmed this by
adding a simulation at that point.

\section{Conclusion and discussion}

We have presented an extension of the method described in~\cite{li08} to the calculation on $\bG^<$ and the current density in NEGF models of nano-transistors. This extension is essential to apply the methods in~\cite{li08} to real engineering problems. This extension assumes that $\bSigma$ has a sparsity pattern similar to $\bA$. It is based on a nested dissection of the mesh and uses a form of Gaussian elimination. In the family of direct methods that provide exact solutions (vs approximate methods), it can be proved that the scaling of the computational cost with problem size using nested dissection is optimal.

We considered the case of a typical finite-difference scheme for a two dimensional device. The sparsity pattern associated with the stencil allows special optimizations, which can reduce the computational cost by a factor of approximately 2. These optimizations are generally applicable to other types of finite-difference stencils, even though this extension was not discussed here. When $\bSigma$ is sparser than $\bA$, for example if it is a diagonal matrix, further cost reductions are possible.

We also presented optimizations for the cases where $\bA$ is symmetric. In this case, the computational cost of many operations can be reduced using Cholesky type factorization. We also proved results regarding properties of positive definiteness for matrices arising in the process, for example $\cU$. Those optimization techniques were derived for computing both $\bG$ and $\bG^<$.

We described a simple approach to reduce the cost of computing the current, which uses an octagonal tiling of the mesh. This allows reducing the number of leaf clusters that are computed while yielding all the off-diagonal entries that are needed.

Finally, regarding the nodes on the physical boundary of the device, we calculated their effect on the computational cost. The cost associated with clusters that share nodes with the physical boundary of the device is typically lower since the size of the boundary set of the cluster is smaller. Consequently, this reduces significantly the computational cost of operating on the clusters near the top of the tree.

Numerical results confirm the computational cost analysis and show the speed-up that can be expected from applying these methods to real-life applications.

We have not discussed the role that certain classes of fast methods like multigrid algorithms can play for our problem. Multigrid algorithms have been used successfully to solve a wide range of linear systems $\bA x=b$ in linear time~\cite{trottenberg2009multigrid,brandt2000general}. In the area of semiconductor process and device simulation, multigrid has been applied, for example, to the solution of various partial differential equations, including Lam\'e equations (linear elasticity), reaction-diffusion (-convection) and drift-diffusion (-convection-reaction) systems~\cite{trottenberg2009multigrid}. In the NEGF area, multigrid has been used in the solution of the Kohn-Sham equations (specifically the associated generalized eigenvalue problem)~\cite{fattebert2000towards,bernholc2008recent,feng2006nonlinear,wijesekera2007efficient}. The computational cost of multigrid when solving a linear system is lower than the FIND algorithm. However, it is not known whether adapting multigrid from solving a linear system of equations $\bA x=b$ to computing in linear time the diagonal of $\bG^r = \bA^{-1}$ and $\bG^r \; \bSigma \; [\bG^r]^\dagger$ is possible.

\section*{\sc Acknowledgement}

We thank the Stanford School of Engineering, the Institute for Computational and Mathematical Engineering, and U.S.\ Army grant
W911NF-07-2-0027 from the Army High Performance Computing Research Center for Agility, Survivability and Informatics at Stanford,
for their support.

\appendix

In the appendix we provide additional technical results that were used in the main body of this paper.

\section{Symmetry preservation}

\label{proof_sympres}

\proof{This property holds for all the leaf clusters by the symmetry of the original matrix $\bA$. For any node in the cluster
tree, if the property holds for its two child clusters $i$ and $j$, then the property holds for their parent node $k$ as well
because the blocks of $\bA$ in~\eqref{eq:major-update-simple} are given by
\begin{subequations}    \label{eq:block-property}
\begin{align}
  \label{eq:block_property1}    
  \bA_{k}(\sB_k, \sB_k) & =
  \begin{pmatrix}
    \bA_i(\sB_k \cap \sB_i, \sB_k \cap \sB_i) &\bA(\sB_k \cap \sB_i, \sB_k \cap \sB_j)                \\
    \bA(\sB_k \cap \sB_j, \sB_k \cap \sB_i) &\bA_j(\sB_k \cap \sB_j, \sB_k \cap \sB_j)
  \end{pmatrix}, \\
  \label{eq:block_property2}
  \bA_{k}(\sS_k, \sS_k) & =
  \begin{pmatrix}
    \bA_i(\sS_k \cap \sB_i, \sS_k \cap \sB_i) &\bA(\sS_k \cap \sB_i, \sS_k \cap \sB_j)
                \\
    \bA(\sS_k \cap \sB_j, \sS_k \cap \sB_i) &\bA_j(\sS_k \cap \sB_j, \sS_k \cap \sB_j)
  \end{pmatrix}, \\
  \label{eq:block_property3}
  \bA_{k}(\sB_k, \sS_k) & =
  \begin{pmatrix}
    \bA_i(\sB_k \cap \sB_i, \sS_k \cap \sB_i) &\b0                               \\
    \b0 &\bA_j(\sB_k \cap \sB_j, \sS_k \cap \sB_j)
  \end{pmatrix}, \\
  \label{eq:block_property4}
  \bA_{k}(\sS_k, \sB_k) & =
  \begin{pmatrix}
    \bA_i(\sS_k \cap \sB_i, \sB_k \cap \sB_i) &\b0                               \\
    \b0 &\bA_j(\sS_k \cap \sB_j, \sB_k \cap \sB_j)
  \end{pmatrix}.
\end{align}
\end{subequations}
}

\section{Positive-definiteness preservation}

\label{proof_posdef}

\proof{Write $n\times n$ symmetric positive definite matrix $\bA$ as $ \begin{pmatrix} a_{11} &z^T    \\ z   &\bar{\bA}_{11}
\end{pmatrix}$ and let $\tilde{\bA}_{11} \lsdef \bar{\bA}_{11} - zz^T/a_{11}$ be the Schur complement of $a_{11}$. Let $\bA^{(0)}
\lsdef \bA$, $\bA^{(1)} \lsdef \tilde{\bA}^{(0)}_{11}$, \ldots, $\bA^{(n-1)} \lsdef \tilde{\bA}^{(n-2)}_{11}$. Note that since
$\bA$ is positive definite, $\bA^{(i)}_{11} \not= 0$ and  $\tilde{\bA}^{(i+1)}$ is well defined for all $i=0, \ldots, n-2$. By
definition, $\tilde{\bA}^{(i+1)}$ are also all symmetric.

Now, we will show that $\bA^{(1)}$, \ldots, $\bA^{(n-1)}$ are all positive definite. Given any $(n-1)\times 1$ matrix $y\not= 0$,
let $x = \begin{pmatrix} -z^Ty/a_{11} \\ y \end{pmatrix} \not= 0$. Since $\bA$ is symmetric positive definite, we have
\begin{eqnarray*}
0   &<  &x^T \bA x    \\
    &= &\begin{pmatrix} -z^T y / a_{11} & y^T \end{pmatrix}
        \begin{pmatrix} a_{11}  &z^T    \\ z    &\bar{\bA}_{11}   \end{pmatrix}
        \begin{pmatrix} -z^Ty/a_{11}    \\  y   \end{pmatrix}   \\
    &= &\begin{pmatrix} 0    &-y^Tzz^T/a_{11} + y^T\bar{\bA}_{11}  \end{pmatrix}   \begin{pmatrix} -z^Ty/a_{11}    \\ y    \end{pmatrix}   \\
    &=  &-y^Tzz^Ty/a_{11} + y^T\bar{\bA}_{11} y    \\
    &=  &y^T(\bar{\bA}_{11} - zz^T/a_{11}) y   \\
    &=  &y^T \tilde{\bA}_{11} y.
\end{eqnarray*}
As a result, $\bA^{(1)}$ is positive definite as well. Repeat this process to show that $\bA^{(2)}$, \ldots, $\bA^{(n-1)}$ are
all positive definite.

Since any principal submatrix of a positive definite matrix is also positive definite~\cite{golub96}, every $\bA(\sS_k, \sS_k)$
in \eqref{eq:major-update} is also positive definite.}

\section{Optimizations that use the sparsity of $\bA$}

The sparsity pattern was described in the main body of the paper. In this section, based on this sparsity pattern and different ways to calculate the inverse of a matrix, we reduce the number of floating point operations.

Recall that we need to calculate the term:
\begin{equation*}
	\bA(\sB, \sS) \; \bA(\sS, \sS)^{-1} \; \bA(\sS, \sB)	
\end{equation*}
which is part of the update of $\cU$:
\begin{equation*}
  \cU = \bA(\sB, \sB) - \bA(\sB, \sS)\bA(\sS, \sS)^{-1}\bA(\sS, \sB)
\end{equation*}

For notational simplicity, in this section, we write
$\bA(\sS, \sS)$ as a block matrix:
$$\blockm{A}{B}{C}{D} = \blockm{\bA(\sS_L, \sS_L)}{\bA(\sS_L, \sS_R)}{\bA(\sS_R, \sS_L)}{\bA(\sS_R, \sS_R)};$$
write $\bA(\sS, \sB)$ as
$$\blockm{W}{X}{Y}{Z} = \blockm{\bA(\sS_L, \sB_L)}{\bA(\sS_L, \sB_R)}{\bA(\sS_R, \sB_L)}{\bA(\sS_R, \sB_R)};$$
and write $\bA(\sB, \sS)$ as
$$\blockm{P}{Q}{R}{S} = \blockm{\bA(\sB_L, \sS_L)}{\bA(\sB_L, \sS_R)}{\bA(\sB_R, \sS_L)}{\bA(\sB_R, \sS_R)}.$$
Since a multiplication is much more expensive than an addition, we ignore the addition with $\bA(\sB, \sB)$.

Before proceeding to the main analysis, we list some facts about the computational cost of various basic matrix computations. We indicate the cost of various operations (shown with $\Rightarrow$) involving the $n\times
n$ full matrices $A$ and $B$ (only multiplications are counted):
\begin{equation*}
	\begin{array}{l@{\quad\Rightarrow\quad}ll}
		A & LU: & n^3/3 \\
		L & L^{-1}: & n^3/6 \\
		L, B & L^{-1}B: & n^3/2 \\
		U, L^{-1}B & A^{-1}B: & n^3/2
	\end{array}
\end{equation*}
Adding them together, computing $A^{-1}$ requires $n^3$ and computing $A^{-1}B$ requires $\frac{4}{3}n^3$. The order of
computation is often important, e.g., if we compute $A^{-1}B = (U^{-1}L^{-1})B$, then the cost is $(n^3 + n^3 =)\ 2n^3$, whereas
computing $U^{-1}(L^{-1}B)$ only requires $\frac{4}{3}n^3$.

We will consider different block forms for the inverse:\footnote{We assume $A$ and $D$ to be nonsingular.}
\begin{subequations}
\begin{eqnarray*}
\begin{pmatrix}
A   &B  \\
C   &D
\end{pmatrix} ^{-1} &= &\begin{pmatrix}
                        \tilde{A}^{-1}              &-A^{-1} B \tilde{D}^{-1}\\
                        -D^{-1} C \tilde{A}^{-1}    &\tilde{D}^{-1}
                        \end{pmatrix} \\
                    &= &\begin{pmatrix}
                        A^{-1} + A^{-1}B\tilde{D}^{-1}CA^{-1}              &-A^{-1} B \tilde{D}^{-1}    \\
                        -\tilde{D}^{-1} C A^{-1}    &\tilde{D}^{-1}
                        \end{pmatrix} \\
                    &= &\begin{pmatrix}
                        A   &B  \\
                        0   &\tilde{D}
                        \end{pmatrix}^{-1}
                        \begin{pmatrix}
                        I   &0  \\
                        CA^{-1}   &I
                        \end{pmatrix}^{-1}
\end{eqnarray*}%
\end{subequations}%
See~\eqref{eq:back-sub-2way}, \eqref{eq:back-sub-1way}, and~\eqref{eq:block-LU-inverse}.
In~\eqref{eq:back-sub-2way}, the calculations for $\tilde{A}^{-1}$ and $\tilde{D}^{-1}$ are independent of each other and can be
done in parallel, so we call this method \emph{parallel inverse}. In~\eqref{eq:back-sub-1way}, we have to calculate
$\tilde{D}^{-1}$ first to calculate the other block of the inverse: $A^{-1} + A^{-1}B\tilde{D}^{-1}CA^{-1}$, so we call this
method \emph{sequential inverse}. In~\eqref{eq:block-LU-inverse}, we perform block LU factorization first and then calculate the
inverse, so we call this method \emph{block LU inverse}.

The first method is based on \eqref{eq:back-sub-2way}, the parallel inverse of the block matrix. It computes $\bA(\sS, \sS)^{-1}$ through the Schur complement of the two
diagonal blocks of $\bA(\sS, \sS)$. We multiply \eqref{eq:back-sub-2way} by $\bA(\sS, \sB)$ on the right. The result of this
multiplication is a $2 \times 2$ block matrix, whose blocks are shown in Table~\ref{block-expressions}. The calculation of each
block requires a set of operations listed in the {\it Operations} columns. The details of how to perform those operations in
order to minimize the computational cost are given in Table~\ref{operation-cost1}.

\begin{table}[htbp]
\centering \caption{Matrix blocks and their corresponding operations} \label{block-expressions}
\begin{tabular}{@{}ccc@{}}    \toprule
\emph{Block}    &\emph{Expression}  &\emph{Operations} \\ \midrule
(1, 1)  &$\tilde{A}^{-1}W - A^{-1} B \tilde{D}^{-1}Y$   &\circlednumber{5} \circlednumber{12}     \\
\addlinespace[2pt]%
(1, 2)  &$\tilde{A}^{-1}X - A^{-1} B \tilde{D}^{-1}Z$   &\circlednumber{9} \circlednumber{8}      \\
\addlinespace[2pt]%
(2, 1)  &$-D^{-1} C \tilde{A}^{-1}W + \tilde{D}^{-1}Y$  &\circlednumber{7} \circlednumber{10}     \\
\addlinespace[2pt]%
(2, 2)  &$D^{-1} C \tilde{A}^{-1}X + \tilde{D}^{-1}Z$   &\circlednumber{11} \circlednumber{6}     \\
\bottomrule
\end{tabular}
\end{table}
Table~\ref{operation-cost1} gives the cost of each operation with different assumptions for the sparsity pattern of $X$ and $Y$.
The last column indicates the dependency between these operations, as some of the operations depend on the results of previous
operations.

In this table and the tables below, to show dependence, the operations are listed in the order of computation. The size of $A$, $B$, $C$, and $D$ is $m\times m$; the size of $W$, $X$, $Y$, and $Z$ is $m\times n$. The different operations are labeled using a number inside a circle. This notation was also used in Table~\ref{block-expressions}.

\begin{table}[htbp]
\centering \caption{The cost of operations and their dependence in the first method. The costs are in flops. The size of $A$,
$B$, $C$, and $D$ is $m\times m$; the size of $W$, $X$, $Y$, and $Z$ is $m\times n$.} \label{operation-cost1}
\begin{tabular}{@{}ccccc@{}}   \toprule
    &\multicolumn{3}{c}{\emph{Type of matrix blocks}}  &    \\
        \cmidrule{2-4}
        \addlinespace[5pt]%
        \emph{Operations}     &all full   &X, Y = 0
        &\begin{minipage}[t]{30pt}\vspace{-12pt}\shortstack{$_{_\text{X, Y = 0;}}$ \\ $_{_{\text{B, C = diag}}}$}\end{minipage}
    &\emph{Dependency}   \\ \hline
         \addlinespace[3pt]%
        \circlednumber{1} \hfill $D\backslash C$     &$4m^3/3$    &$4m^3/3$    &$m^3$    &n/a    \\
        \addlinespace[2pt]%
        \circlednumber{2} \hfill $A\backslash B$     &$4m^3/3$    &$4m^3/3$    &$m^3$          &n/a    \\      
        \addlinespace[2pt]%
        \circlednumber{3} \hfill $\tilde{A} = A - BD^{-1}C$      &$m^3$      &$m^3$            &0    &\circlednumber{1}    \\      
        \addlinespace[2pt]%
        \circlednumber{4} \hfill $\tilde{D} = D - CA^{-1}B$      &$m^3$      &$m^3$            &0    &\circlednumber{2}    \\      
        \addlinespace[2pt]%
        \circlednumber{5} \hfill $\tilde{A}\backslash W$    &$\frac{m^3}{3} + m^2n$     &$\frac{m^3}{3} + m^2n$
            &$\frac{m^3}{3} + m^2n$        &\circlednumber{3}    \\      
        \addlinespace[2pt]%
        \circlednumber{6} \hfill $\tilde{D}\backslash Z$    &$\frac{m^3}{3} + m^2n$    &$\frac{m^3}{3} + m^2n$
            &$\frac{m^3}{3} + m^2n$        &\circlednumber{4}    \\      
        \addlinespace[2pt]%
        \circlednumber{7} \hfill $-D^{-1}C\tilde{A}^{-1}W$      &$m^2n$      &$m^2n$      &$m^2n$   &\circlednumber{1}\circlednumber{5}\\   
        \addlinespace[2pt]%
        \circlednumber{8} \hfill $-A^{-1}B\tilde{D}^{-1}Z$      &$m^2n$      &$m^2n$      &$m^2n$   &\circlednumber{2}\circlednumber{6}\\  
        \addlinespace[2pt]%
        \circlednumber{9} \hfill $\tilde{A}^{-1}X$  &$m^2n$      &0      &0          &\circlednumber{5}    \\      
        \addlinespace[2pt]%
        \circlednumber{10} \hfill $\tilde{D}^{-1}Y$ &$m^2n$      &0      &0         &\circlednumber{6}    \\      
        \addlinespace[2pt]%
        \circlednumber{11} \hfill $-D^{-1}C\tilde{A}^{-1}X$     &$m^2n$      &0      &0  &\circlednumber{1}\circlednumber{9}\\      
        \addlinespace[2pt]%
        \circlednumber{12} \hfill $-A^{-1}B\tilde{D}^{-1}Y$     &$m^2n$      &0      &0  &\circlednumber{2}\circlednumber{10}\,\\      
        \addlinespace[2pt]%
        \hline
        \addlinespace[2pt]%
        \textbf{Total}                 &$\frac{16m^3}{3}+8m^2n$    &$\frac{16m^3}{3}+4m^2n$    &$\frac{8m^3}{3}+4m^2n$      &n/a \\
        \addlinespace[2pt]%
        \bottomrule
    \end{tabular}
\end{table}
Since the clusters are typically equally split, we can let $m=|\sS_L| = |\sS_R| = |\sS|/2$ and $n = |\sB_L| = |\sB_R| = |\sB|/2$.
Then, the size of matrices $A, B, C$, and $D$ is $m\times m$; the size of matrices $W, X, Y$, and $Z$ is $m\times n$, and the
size of matrices $P, Q, R$, and $S$ is $m\times n$. For example, for a merge from two $(a\times a)$-clusters to one $(a\times
2a)$-cluster, we have $m=a$ and $n=3a$.

Table~\ref{operation-cost1} gives the cost of computing $\bA(\sS, \sS)^{-1}\bA(\sS, \sB)$ in \eqref{eq:major-update-simple}.
Consider now the multiplication by $\bA(\sB, \sS)$ from the left: it requires $4m^2n$ operations. The total computation cost for
\eqref{eq:major-update-simple} with the above optimization is therefore $\frac{8}{3}m^3 + 4m^2n + 4mn^2$. In contrast, the cost
without optimization is $8mn^2 + 8m^2n + \frac{8}{3}m^3$ since $s=2m$ and $b=2n$.

Taking two upward cases as examples, $(m, n) = (a, 3a)$ and $(m, n) = (2a, 4a)$, the reduction is:
\begin{equation*}
	\frac{296}{3}a^3 \rightarrow \frac{152}{3}a^3 \quad
	\text{ and } \quad
	\frac{1216}{3}a^3 \rightarrow \frac{640}{3}a^3.
\end{equation*}
The cost with optimization lies somewhere between 51\% to 53\% of the original cost. Similarly, for the downward pass, in the two
typical cases $(m, n) = (3a, 2a)$ and $(m, n) = (4a, 3a)$, we have the reduction
\begin{equation*}
	312a^3 \rightarrow 192a^3 \; \text{(62\%)}\quad
	\text{ and } \quad
	\frac{2528}{3} \rightarrow \frac{1520}{3}a^3 \; \text{(60\%)}.
\end{equation*}
In the downward pass, $|\sB_L| \neq |\sB_R|$, but in terms of the computational cost estimate, we can assume that they are both $s/2$.

The second method is based on \eqref{eq:back-sub-1way}, the sequential inverse of the block matrix. Similar to the first method, it does not compute $\bA(\sS, \sS)^{-1}$
explicitly. Instead, it computes $\bA(\sS, \sS)^{-1}\bA(\sS, \sB)$ as
\begin{equation}
{\footnotesize\begin{pmatrix}
A   &B  \\
C   &D
\end{pmatrix} ^{-1}
\begin{pmatrix}
W  &0   \\
0   &Z
\end{pmatrix}
=
\begin{pmatrix}
A^{-1}W + A^{-1}B\tilde{D}^{-1}CA^{-1}W     &-A^{-1}B\tilde{D}^{-1}Z    \\
-\tilde{D}^{-1}CA^{-1}W                     &\tilde{D}^{-1}Z
\end{pmatrix}},
\label{eq:block-inverse-multi}
\end{equation}
where $\tilde{D} = D-CA^{-1}B$, $B$ and $C$ are diagonal, and $X=Y=0$. Table~\ref{operation-cost2} shows the required operations
with the total cost $\frac{4}{3}m^3 + 5m^2n + 4mn^2$.

\begin{table}[htbp]
\centering \caption{Operation costs in flops and their dependence in the sequential inverse method.} \label{operation-cost2}
\begin{tabular}{@{}cccc@{}}    \toprule
    \emph{Operation}     &\emph{Cost}     &\emph{Dependency}   \\ \midrule
    \addlinespace[3pt]%
    \circlednumber{1} \hfill $A^{-1}$       &$m^3$              &n/a        \\      
    \addlinespace[2pt]%
    \circlednumber{2} \hfill $\tilde{D}$    &0                     &\circlednumber{1}     \\  
    \addlinespace[2pt]%
    \circlednumber{3} \hfill $\tilde{D}^{-1}Z$           &$\frac{m^3}{3}+m^2n$      &\circlednumber{2}    \\  
    \addlinespace[2pt]%
    \circlednumber{4} \hfill $-(A^{-1}B)(\tilde{D}^{-1}Z)$   &$m^2n$                    &\circlednumber{1}\circlednumber{3} \\  
    \addlinespace[2pt]%
    \circlednumber{5} \hfill $A^{-1}W$                  &$m^2n$                    &\circlednumber{1}     \\      
    \addlinespace[2pt]%
    \circlednumber{6} \hfill $-\tilde{D}^{-1}(CA^{-1}W)$    &$m^2n$                    &\circlednumber{3}\circlednumber{5} \\  
    \addlinespace[2pt]%
    \circlednumber{7} \hfill $(A^{-1}W) + (A^{-1}B)(\tilde{D}^{-1}CA^{-1}W)$ &$m^2n$  &\circlednumber{1}\circlednumber{5}\circlednumber{6}    \\  
    \addlinespace[2pt]%
    \circlednumber{8} \hfill $\bA(\sS, \sB)\bA(\sS, \sS)^{-1}\bA(\sS, \sB)$  &$4mn^2$      &\circlednumber{7}     \\
    \addlinespace[2pt]%
    \hline
    \addlinespace[2pt]%
    \textbf{Total}  &$\frac{4}{3}m^3 + 5m^2n + 4mn^2$       &n/a        \\      \bottomrule
\end{tabular}
\end{table}

Compared to the previous method, the cost here is smaller if $m>\frac{3}{4}n$. Therefore, this method is better than the first
method for the downward pass where typically $(m, n)$ is $(3a, 2a)$ and $(4a, 2a)$. In such cases, the costs are $174a^3$ (56\%)
and $\frac{1408}{3}a^3$ (56\%), respectively.

The third method is based on \eqref{eq:block-LU-inverse}, the block inverse of the block matrix. The two factors there will multiply $\bA(\sB, \sS)$ and $\bA(\sS, \sB)$
separately:
\begin{eqnarray}
\lefteqn{\bA(\sB, \sS) \; \bA(\sS, \sS)^{-1} \; \bA(\sS, \sB) } \notag \\
    &= &\left[ \blockm{P}{0}{0}{S} \blockm{A}{B}{0}{\tilde{D}}^{-1} \right]
        \left[\blockm{I}{0}{CA^{-1}}{I}^{-1} \blockm{W}{0}{0}{Z} \right]    \notag  \\
    &= &\left[ \blockm{P}{0}{0}{S} \blockm{A^{-1}} {-A^{-1}B\tilde{D}^{-1}} {0} {\tilde{D}^{-1}} \right]
        \left[\blockm{I}{0}{-CA^{-1}}{I} \blockm{W}{0}{0}{Z} \right]   \notag  \\
    &= &\left[ \blockm{PA^{-1}}{-PA^{-1}B\tilde{D}^{-1}}{0}{S\tilde{D}^{-1}} \right]
        \left[ \blockm{W}{0}{-CA^{-1}W}{Z} \right]  \label{eq:block-LU-multi}
\end{eqnarray}
The cost of \eqref{eq:block-LU-multi} is $\frac{4}{3}m^3 + 4m^2n + 5mn^2$, which comes from the multiplication of a block upper triangular matrix and a block lower triangular matrix. Table \ref{operation-cost3} shows the order and details of the operations. $S\tilde{D}^{-1}$ is computed by first LU factorizing $\tilde{D}$ and then solving for $S\tilde{D}^{-1}$. The LU form of $\tilde{D}$ will be used when computing $-(PA^{-1}B)\tilde{D}^{-1}$ as well.

\begin{table}[htbp]
\centering \caption{Operation costs in flops and their dependence in the block LU inverse method.} \label{operation-cost3}
\begin{tabular}{@{}cccc@{}} \toprule
     \emph{Operation}     &\emph{Cost}     &\emph{Dependency}   \\ \midrule
     \addlinespace[3pt]%
     \circlednumber{1} \hfill $A^{-1}$       &$m^3$              &n/a        \\
     \addlinespace[2pt]%
     \circlednumber{2} \hfill $\tilde{D}$    &0                     &\circlednumber{1}     \\
     \addlinespace[2pt]%
     \circlednumber{3} \hfill $\tilde{D}^{-1}Z$           &$\frac{m^3}{3}+m^2n$      &\circlednumber{2}    \\
     \addlinespace[2pt]%
     \circlednumber{4} \hfill $-(A^{-1}B)(\tilde{D}^{-1}Z)$   &$m^2n$                    &\circlednumber{1}\circlednumber{3} \\
     \addlinespace[2pt]%
     \circlednumber{5} \hfill $A^{-1}W$                  &$m^2n$                    &\circlednumber{1}     \\
     \addlinespace[2pt]%
     \circlednumber{6} \hfill $-\tilde{D}^{-1}(CA^{-1}W)$    &$m^2n$                    &\circlednumber{3}\circlednumber{5} \\
     \addlinespace[2pt]%
     \circlednumber{7} \hfill $(A^{-1}W) + (A^{-1}B)(\tilde{D}^{-1}CA^{-1}W)$ &$m^2n$  &\circlednumber{1}\circlednumber{5}\circlednumber{6}    \\
     \addlinespace[2pt]%
     \circlednumber{8} \hfill $\bA(\sS, \sB)\bA(\sS, \sS)^{-1}\bA(\sS, \sB)$  &$4mn^2$      &\circlednumber{7}     \\
     \addlinespace[2pt]%
    \hline
     \addlinespace[2pt]%
     \textbf{Total}  &$\frac{4}{3}m^3 + 5m^2n + 4mn^2$        &n/a        \\      \bottomrule
\end{tabular}
\end{table}

Finally, our fourth method is more straightforward. We call it \emph{na\"{i}ve LU method}. It considers $\bA(\sB, \sS)$ and $\bA(\sS, \sB)$ as block diagonal matrices but considers $\bA(\sS, \sS)$ as a full matrix without exploiting its sparsity. This method gives cost $\frac{8}{3}m^3 + \frac{13}{2}m^2n + 4mn^2$. Table~\ref{operation-cost4} shows the order and details of the operations. The cost for computing the second column of $L^{-1}\bA(\sS, \sB)$ is reduced since $\bA(\sS, \sB)$ is block diagonal.

\begin{table}[htbp]
\centering \caption{Operation cost in flops and their dependence in the na\"{i}ve LU inverse method.}\label{operation-cost4}
\begin{tabular}{@{}cc@{}}    \toprule
     \emph{Operation}     &\emph{Cost}     \\ \midrule
     \addlinespace[2pt]%
     LU for $\bA(\sS, \sS)$     &$\frac{8}{3}m^3$   \\
     \addlinespace[2pt]%
     $L^{-1}\bA(\sS, \sB)$: 1st block column      &$2m^2n$    \\
     \addlinespace[2pt]%
     $L^{-1}\bA(\sS, \sB)$: 2nd block column      &$\frac{1}{2}m^2n$    \\
     \addlinespace[2pt]%
     $U^{-1}[L^{-1}\bA(\sS, \sB)]$           &$4m^2n$    \\
     \addlinespace[2pt]%
     $\bA(\sS, \sB)\bA(\sS, \sS)^{-1}\bA(\sS, \sB)$  &$4mn^2$  \\
     \addlinespace[2pt]%
    \hline
     \addlinespace[2pt]%
     \textbf{Total}  &$\frac{8}{3}m^3 + \frac{13}{2}m^2n + 4mn^2$     \\  \bottomrule
\end{tabular}
\end{table}

The summary of all the optimization techniques for sparsity and the corresponding costs is given in the main body of this paper, in Table~\ref{operation-cost-summary}.

Finally, Table~\ref{NullBoundaryCost} shows the estimated level-by-level costs with the physical boundary of the device taken into consideration (see Section ``Nodes on the physical boundary of the device''). This cost does not assume any optimization. It simply accounts for the effect of the boundary of the physical device, which leads to a computational cost reduction.

The last two rows for each pass are the sum of the cost of the rest of the small clusters. This estimate does not consider the two end blocks corresponding to the device contacts as in real applications~\cite{li08,svizhenko02}. Accounting for the two end blocks, if we do not have any optimization, then the grand total is: $450N^3 + \frac{16}{3}N^3 + 2N^3 \approx 457 N^3$. As a comparison, the cost using RGF with symmetry taken into consideration is $3.5N^4$. Since the cross-point of the two cost curves is at $N$ equal to the ratio between the constant factors in the two methods, it is expected to happen at $N\thicksim130$. This is consistent with the results from our early simulations (between $N=100$ and $N=150$)~\cite{li08}.

\begin{table}[htp]
\begin{center}
\caption {Computation cost estimate for different type of clusters from an $N\times N$ mesh. $\Ss$ is the size of the private
inner node set and $\Sb$ is the size of the boundary set. The size is in unit of $N$ and the cost is in unit of $N^3$. By
convention, the root (the whole mesh) of the cluster tree is level 0. The cost column gives the total cost for the given subset
of clusters. For example, on row one, under upward pass, we have 2 clusters at level 1 with $\Ss=1 \, N$ and $\Sb=1 \, N$ and the
cost for operating on these clusters is 4.666\,$N^3$.} \label{NullBoundaryCost}
\begin{tabular}{@{}cccccccccc@{}}    \toprule
\multicolumn{5}{c}{\emph{upward pass}}    &\multicolumn{5}{c}{\emph{downward pass}}  \\
\cmidrule(r){1-5}   \cmidrule(l){6-10}
\multicolumn{3}{c}{\emph{size per cluster \& total cost}}
&\multicolumn{2}{c}{\emph{clusters}}
&\multicolumn{3}{c}{\emph{size \& cost per cluster}}
&\multicolumn{2}{c}{\emph{clusters}}   \\
\cmidrule(r){1-3} \cmidrule(r){4-5} \cmidrule(r){6-8} \cmidrule(r){9-10}
$\Ss$ &$\Sb$  &\emph{cost} &\emph{level} &\emph{number} &$\Ss$ &$\Sb$  &\emph{cost} &\emph{level} &\emph{number} \\
\midrule
1   &1  &4.666  &1  &2
&0   &1  &0      &1  &2      \\
1   &1  &9.332  &2  &4
&1       &1      &9.332      &2  &4     \\
1/2   &3/4  &2.040  &3  &4
&3/2      &3/4  &14.624   &3  &4      \\
1/2     &5/4  &4.540  &3  &4
&1/2       &5/4     &4.540   &3  &4      \\
1/2   &1  &3.168    &4  &4
&1       &1/2      &4.332    &4  &4     \\
1/2   &3/4  &4.080  &4  &8
&1/2       &3/4     &2.040  &4  &4      \\
1/2     &1/2  &1.168  &4  &4
&3/2      &3/4      &14.624   &4  &4      \\
1/4   &3/8  &0.2552  &5  &4
&1       &1      &9.332      &4  &4     \\
1/4   &1/2  &0.3958  &5   &4
&1       &3/4   &52.672    &5  &32     \\
1/4     &5/8  &1.703  &5  &12
&3/4       &1/2      &38.976  &6  &64     \\
1/4   &3/4  &2.312  &5  &12
&       &      &$<$91.648    &7...  &     \\
1/4   &1/2  &6.336  &6  &64     \\
1/8   &3/8  &3.072  &7  &128    \\
& & $<$9.408 &8...  & \\
\midrule
\multicolumn{2}{l}{\textbf{Subtotal}}      &52.470  &   &
&\multicolumn{2}{l}{\textbf{Subtotal}}      &242.120  &   & \\
\multicolumn{2}{l}{\textbf{Total:}}   &\textbf{294.6}   &   &     \\
\bottomrule
\end{tabular}
\end{center}
\end{table}

\noindent {\it Summary.}\ \ The optimization for sparsity leads to $40\%$ cost reduction and that for symmetry leads to $50\%$ cost reduction. This results in a $70\%$ cost reduction. In other words, the new cost is about $1/3$ of the original cost. This is reflected in the cross-point of the performance curves using these two methods. Originally, the cross-point was around $N\thicksim130$, but now it can be as low as $N=40$. This shows the performance improvement resulting from our optimizations. The optimization for sparsity does not reduce the memory cost. The optimization for symmetry does reduce it (approximately by half).

\section{Additional proofs} \label{tech-supplement}

In this section, we present some additional technical results and proofs required to show correctness of our algorithm. All the relevant properties of the mesh node subsets for our algorithm are listed, along with their proofs for the basic FIND algorithm in Li et al.~\cite{li08}. We will then use these properties to prove two theorems and four corollaries used for computing $\bG^<$. Some of the properties are not directly used in the proofs but are still listed to help understand the algorithm.

Property~\ref{prop:case-by-case} is fairly simple but will be used again and again in proving theorems, corollaries, and other
properties.

\begin{property}\label{prop:case-by-case}
By definition of $\sT_r^+$, one of the following relations must hold: $\sC_i \subset \sC_j$, $\sC_i \supset \sC_j$, or $\sC_i
\cap \sC_j = \emptyset$.
\end{property}

Property~\ref{prop:private-inner-nodes-union} shows an alternative way of looking at the inner mesh nodes. It helps understand
Property~\ref{prop:bb=bu}.

\begin{property}\label{prop:private-inner-nodes-union}
$\sI_i = \cup_{\sC_j\subseteq \sC_i} \sS_j$, where $\sC_i$ and $\sC_j \in \sT_r^+$.
\end{property}

\begin{property}\label{prop:bb=bu}
If $\sC_i$ and $\sC_j$ are the two children of $\sC_k$, then $\sS_k\cup \sB_k = \sB_i \cup \sB_j$ and $\sS_k = (\sB_i\cup
\sB_j)\backslash \sB_k$.
\end{property}

\begin{property}\label{prop:disjoint}
For any given augmented tree $\sT_r^+$ and all $\sC_i \in \sT_r^+$, $\sS_i$, $\sB_{-r}$, and $\sC_r$ are all disjoint and
$\sM=(\cup_{\sC_i \in \sT_r^+} \sS_i) \cup \sB_{-r} \cup \sC_r $.
\end{property}

Property~\ref{prop:disjoint} shows that the whole mesh can be partitioned into subsets $\sS_i$, $\sB_{-r}$, and $\sC_r$. See Sections ``Brief description of the FIND algorithm'' and ``Formal description of the FIND algorithm'' for notations and figures.

Below we list properties of $\sS_i$ for specific orderings.

\begin{property}\label{prop:dichotomy}
If $\sS_i<\sS_j$, then $\sC_j \not\subset \sC_i$, which implies either $\sC_i\subset \sC_j$ or $\sC_i\cap \sC_j = \emptyset$.
\end{property}

This property is straightforward from the definition of $\sS_i<\sS_j$ and Property~\ref{prop:case-by-case}.

Properties~\ref{prop:prop6} and~\ref{prop:no-sandwich} below are related to the elimination process. They were used in the proof of the following theorem (see main body of paper, theorem~\ref{thm:FIND-major}):
\begin{thm}
For any target clusters $r$ and $s$ such that $\sC_g\in \sT_{r}^+$ and $\sC_g\in \sT_{s}^+$, we have
\begin{equation*}
  \bA_{r, g+}(\sB_g, \sB_g) = \bA_{s, g+}(\sB_g, \sB_g).
\end{equation*}
\end{thm}
They will also be used in the proof of Theorem~\ref{thm:individual-update} and Theorem~\ref{thm:synchronized-update}.

\begin{property}\label{prop:prop6}
For any $k, u$ such that $\sC_k, \sC_u \in \sT_r^+$, if $\sC_u\not\subset\sC_k$, then $\sB_u \cap \sI_k = \emptyset$.
\end{property}

\proof{By Property~\ref{prop:case-by-case}, we have $C_u \not\subset C_k$. So for all $j$ such that $C_j \subseteq C_k$, we have
$j\not = u$ and thus $S_j \cap S_u = \emptyset$ by property 3. By property 2, $I_k = \cup_{C_j\subseteq C_k}S_j$, so we have
$I_k\cap S_u = \emptyset$.}

\begin{property}\label{prop:no-sandwich}
If $\sC_j$ is a child of $\sC_k$, then for any $\sC_u$ such that $\sS_j<\sS_u<\sS_k$, we have $\sC_j\cap \sC_u = \emptyset$ and
thus $\sB_j \cap \sB_u = \emptyset$.
\end{property}

This is because the $\sC_u$ can be neither a descendant of $\sC_j$ nor an ancestor of $\sC_k$.

In Theorem~\ref{thm:individual-update} and its proof below, we follow the same notation used in Li et al.~\cite{li08}. In
particular, we write $\sS_g$ more formally as $\sS_{i_j}$, with $j = 1, 2, \ldots$, etc., where $\bA_{i_1} = \bA$, $\bSigma_{i_1}
= \bSigma$, and $\sS_{i_j} < \sS_{i_{j'}}$ iff $j<j'$. We also denote $S_{i_{j+1}}$ as $S_{g+}$ when $g=i_j$. Because of
Property~\ref{prop:disjoint}, the whole mesh $\sM$ is partitioned into a sequence of sets $S_{i_j}$, $\sB_{-r}$, and $\sC_r$,
$j=1, 2, \ldots$, and the order of $\bA$ and $\bSigma$ is determined by this sequence. For notation convenience, if $i$ is the
index of the last set in the sequence of $\sS_{i_j}$, which is always $-r$ for $\sT_r^+$, then $\sS_{i+}$ stands for $\sB_{-r}$.

\begin{thm} \label{thm:individual-update}
If we perform Gaussian elimination as described in the Section ``Formal description of the FIND algorithm'' and update $\bSigma_g$ accordingly based on
$\bSigma_{g+} = \bL^{-1}_g \bSigma_g \bL^{-\dagger}_g$, then we have:
\begin{enumerate}
\item $\forall \sS_h\geq \sS_g, \bSigma_g(\sS_h, \sS_{>h}\backslash \sB_h) = \bSigma_g(\sS_{>h}\backslash \sB_h , \sS_h) =
    0$;
\item
    \begin{enumerate}
    \item $\bSigma_{g+}(\sS_{\leq g}, \sS_{\leq g}) = \bSigma_g(\sS_{\leq g}, \sS_{\leq g})$;
    \item $\bSigma_{g+}(\sM, \sS_{>g}\backslash \sB_g) = \bSigma_g(\sM, \sS_{>g}\backslash \sB_g)$;
    \item $\bSigma_{g+}(\sS_{>g}\backslash \sB_g, \sM) = \bSigma_g(\sS_{>g}\backslash \sB_g, \sM)$;
    \end{enumerate}
\item $\bSigma_{g+}(\sB_g, \sB_g) = \ \bSigma_g(\sB_g, \sB_g) - \cL_g\bSigma_g(\sS_g, \sB_g)
		- \bSigma_g(\sB_g, \sS_g)\cL_g^{\dagger} + \cL_g\bSigma_g(\sS_g, \sS_g)\cL_g^{\dagger}$.
\end{enumerate}
\end{thm}

In other word, this theorem describes how the multiplications by $\bL_g^{-1}$ and $\bL_g^{-\dagger}$ affect $\bSigma_g$.
Specifically, these multiplications only affect $\bSigma_g(\sB_g, \sB_g)$, which is given by (III), and $\bSigma_g(\sS_{\leq g},
\sB_g)$ and $\bSigma_g(\sB_g, \sS_{\leq g})$, which are not given here but can be ignored. The multiplication preserves the
sparsity pattern specified by (I) above, which is illustrated by the following submatrix for any $\sS_h >
\sS_g$:\footnote{Similar to the notation $\bA_g(\bullet, \sB_g)$, the entries of the blocks in
this submatrix usually do \emph{not} stay together (except the entries of $\bSigma_g(\sS_h, \sS_h)$, which always stay together).
We write them here like blocks to better illustrate the sparsity pattern of $\bSigma_g$.}
\begin{equation*}
\begin{array}{cccc}
  &\sS_{h}     &\sB_{h}       &\sS_{>h}\backslash \sB_{h}  \\
  \sS_{h}     &\times  &\times    &\b0                      \\
  \sB_{h}     &\times  &\times    &\times                \\
  \sS_{>h}\backslash B_{h} &\b0 &\times   &\times
\end{array}
\end{equation*}

\proof{Recall that $-\cL_g=-\bL_g(\sB_g, \sS_g)$ is the only non-zero off-diagonal block of $\bL_g^{-1}$.
As a result, all entries in $\bSigma_g(\bullet, \sS_{<g})$ and $\bSigma_g(\sS_{<g}, \bullet)$ are irrelevant to the
multiplication and (I) implies (II) and (III). So it suffices to prove (I) by mathematical induction.

For $g = i_1$, (I) holds because no nodes in $\sS_h$ and $\sS_{>h}\backslash \sB_h$ are connected to each other and in the
original matrix $\bSigma$, an entry is non-zero iff the two nodes associated with the column and the array of the entry connect
to each other.
If (I) holds for $g=k$, then by Property~\ref{prop:case-by-case} we have either $\sC_k \subset \sC_{k+}$ or $\sC_k \cap \sC_{k+}
= \emptyset$. So,
\begin{lsitemizei}
\item if $\sC_k \subset \sC_{k+}$, consider $u$ such that $\sS_{k+} < \sS_u$. By Property~\ref{prop:prop6}, we have $\sI_{k+}
    \cap \sS_u = \emptyset$. Since $\sC_{k+}=\sI_{k+}\cup \sB_{k+}$, we have $(\sS_u\backslash \sB_{k+})\cap \sC_{k+} = \emptyset$. So we have
  $\sB_k \cap (\sS_u\backslash \sB_{k+}) \subseteq (\sS_u\backslash \sB_{k+})\cap \sC_k \subseteq (\sS_u\backslash \sB_{k+})\cap
  \sC_{k+} = \emptyset \Rightarrow \sB_k \cap (\sS_{>k+}\backslash \sB_{k+}) = \emptyset$;
\item if $\sC_k \cap \sC_{k+} = \emptyset$, then $\sB_k\subset \sC_k$ and $\sS_{k+}\subset \sC_{k+} \Rightarrow \sB_k\cap
    \sS_{k+} = \emptyset$.
\end{lsitemizei}
So in both cases, we have $(\sB_k, \sB_k) \cap (\sS_{k+}, \sS_{>k+}\backslash \sB_{k+}) = \emptyset$. In addition, by
Property~\ref{prop:disjoint}, we have $\sS_k \cap \sS_{k+} = \emptyset$. Since all (I), (II), and (III) hold for $g=k$, we have
$\bSigma_{k+}(\sS_{k+}, \sS_{>k+}\backslash \sB_{k+}) = \bSigma_{k+}(\sS_{>k+}\backslash \sB_{k+}, \sS_{k+}) = \b0$.

Since the above argument is valid for all $h\geq k+$, we have $\forall h \geq k+$, $\bSigma_{k+}(\sS_h, \sS_{>h}\backslash
\sB_h)= \bSigma_{k+}(\sS_{>h}\backslash \sB_h, \sS_h) = \b0$, i.e., (I) holds for $g=k+$ as well. By induction, we have that (I)
holds for all $g$ such that $\sC_g \in \sT_r^+$.}

\begin{cor} \label{cor:unchanged-entries} 
If $\sC_i$ and $\sC_j$ are the two children of $\sC_k$, then $\bSigma_k(\sB_i, \sB_j) = \bSigma(\sB_i, \sB_j)$ and $\bSigma_k(\sB_j, \sB_i) = \bSigma(\sB_j, \sB_i)$.
\end{cor}

In words, the corollary states that when we compute $\bL_k^{-1}\bSigma_k\bL_k^{-\dagger}$, we can take some of the entries in
$\bSigma_k$ directly from the original $\bSigma$.   \\

\proof{By Theorem~\ref{thm:individual-update}, the multiplications by $\bL_u^{-1}$ and $\bL_u^{-\dagger}$ only change the entries
of $\bSigma_u$ in $(\sB_u, \sB_u)$, $(\sS_{\leq u}, \sB_u)$, and $(\sB_u, \sS_{\leq u})$. So it suffices to show that for every
$\sS_u < \sS_k$, we have
\begin{equation}
(\sB_u, \sB_u) \cap (\sB_i, \sB_j) = \emptyset,  \label{eq:bb-identity}
\end{equation}
\begin{equation}
(\sS_{\leq u}, \sB_u) \cap (\sB_i, \sB_j) = (\sB_u, \sS_{\leq u})\cap (\sB_i, \sB_j) = \emptyset,        \label{eq:sb-identity}
\end{equation}
and similar identities with $i$ and $j$ interchanged.

To show~\eqref{eq:bb-identity}, consider the following three cases:
\begin{lsitemizei}
\item $\sC_u \cap \sC_k = \emptyset \Rightarrow \sC_u \cap \sC_i = \sC_u \cap \sC_j = \emptyset
    \Rightarrow \sB_u \cap \sB_i = \sB_u \cap \sB_j = \emptyset$;
\item $\sC_u \subset \sC_i \Rightarrow \sC_u \cap \sC_j = \emptyset \Rightarrow \sB_u \cap \sB_j = \emptyset$;
\item $\sC_u \subset \sC_j \Rightarrow \sC_u \cap \sC_i = \emptyset \Rightarrow \sB_u \cap \sB_i = \emptyset$.
\end{lsitemizei}
Because $\forall \sS_u < \sS_k$, by Property~\ref{prop:dichotomy} (with some extension), one of the above cases must hold. So we
have proved~\eqref{eq:bb-identity}.

To show~\eqref{eq:sb-identity}, consider $\sS_v < \sS_k$. Because
$$\sS_v < \sS_k \Rightarrow \sC_k\not\subset\sC_v \Rightarrow
\sC_i\not\subset\sC_v$$ by Property~\ref{prop:dichotomy}, we have either
$$\sC_v\cap\sC_i = \emptyset \Rightarrow \sS_v\cap\sB_i = \emptyset$$
or
$$\sC_v\subset\sC_i \mbox{ (by Property~\ref{prop:private-inner-nodes-union})} \Rightarrow \sS_v\subset\sI_i \Rightarrow \sS_v\cap
\sB_i = \emptyset.$$ We can also conclude $\sS_v\cap\sB_j = \emptyset$ by a similar argument, so we have proved \eqref{eq:sb-identity}.
}

\begin{cor} \label{cor:changed-entries}
If $\sC_i$ is a child of $\sC_k$, then $\bSigma_k(\sB_i, \sB_i) = \bSigma_{i+}(\sB_i, \sB_i)$.
\end{cor}
\proof{Consider $u$ such that $\sS_i < \sS_u < \sS_k$. By Property~\ref{prop:no-sandwich}, we have $\sB_i \cap \sB_u =
\emptyset$. By Theorem~\ref{thm:individual-update}, the multiplications by $\bL_u^{-1}$ and $\bL_u^{-\dagger}$ only change the
entries of $\bSigma_u$ in $(\sB_u, \sB_u)$, $(\sS_{\leq u}, \sB_u)$, and $(\sB_u, \sS_{\leq u})$ and none of them will change the
entries in $(\sB_i, \sB_i)$, so the corollary is proved.}

\begin{cor} \label{cor:leaf-node} If $\sC_i$ is a leaf node in $\sT_r^+$, then
$\bSigma_i(\sC_i, \sC_i) = \bSigma(\sC_i, \sC_i)$.
\end{cor}
\proof{Consider $\sS_u < \sS_i$. By Property~\ref{prop:case-by-case}, either $\sC_u \subset \sC_i$ or $\sC_u \cap \sC_i =
\emptyset$. Since $\sC_i$ is a leaf node, there is no $u$ such that $\sC_u \subset \sC_i$. So we have $\sC_u \cap \sC_i =
\emptyset$ and thus $\sB_u \cap \sC_i = \emptyset$. By theorem \ref{thm:individual-update}, the corollary is proved.}

\begin{thm} \label{thm:synchronized-update} For any $r$ and $s$ such that $\sC_i\in \sT_{r}^+$ a
nd $\sC_i\in \sT_{s}^+$, we have:
\begin{equation}    \label{eq:synchronized-update}
  \bSigma_{r, i}(\sS_i\cup \sB_i, \sS_i\cup \sB_i) = \bSigma_{s, i}(\sS_i\cup \sB_i, \sS_i\cup
  \sB_i)
\end{equation}
\end{thm}
\proof{If $\sC_i$ is a leaf node, then by Corollary~\ref{cor:leaf-node}, we have $\bSigma_{r, i}(\sS_i\cup \sB_i, \sS_i\cup
\sB_i) = \bSigma_{r, i}(\sC_i, \sC_i) = \bSigma_r(\sC_i, \sC_i) = \bSigma_s(\sC_i, \sC_i) = \bSigma_{s, i}(\sC_i, \sC_i) =
\bSigma_{s, i}(\sS_i\cup \sB_i, \sS_i\cup \sB_i)$

If~\eqref{eq:synchronized-update} holds for $i$ and $j$ such that $\sC_i$ and $\sC_j$ are the two children of $\sC_k$, then
\begin{lsitemizei}
\item by Theorem~\ref{thm:individual-update}, we have
    \begin{lsitemizeii}
    \item $\bSigma_{r, i+}(\sB_i, \sB_i) = \bSigma_{s, i+}(\sB_i, \sB_i)$ and
    \item $\bSigma_{r, j+}(\sB_j, \sB_j) = \bSigma_{s, j+}(\sB_j, \sB_j)$;
    \end{lsitemizeii}
\item by Corollary~\ref{cor:unchanged-entries}, we have
    \begin{lsitemizeii}
    \item $\bSigma_{r, k}(\sB_i, \sB_j) = \bSigma_r(\sB_i, \sB_j) = \bSigma_s(\sB_i, \sB_j) = \bSigma_{s, k}(\sB_i, \sB_j)$ and
    \item $\bSigma_{r, k}(\sB_j, \sB_i) = \bSigma_r(\sB_j, \sB_i) = \bSigma_s(\sB_j, \sB_i) = \bSigma_{s, k}(\sB_j, \sB_i)$.
    \end{lsitemizeii}
\item by Corollary~\ref{cor:changed-entries}, we have
    \begin{lsitemizeii}
    \item $\bSigma_{r, k}(\sB_i, \sB_i) = \bSigma_{r, i+}(\sB_i, \sB_i) = \bSigma_{s, i+}(\sB_i, \sB_i) = \bSigma_{s, k}(\sB_i, \sB_i)$ and
    \item $\bSigma_{r, k}(\sB_j, \sB_j) = \bSigma_{r, j+}(\sB_j, \sB_j) = \bSigma_{s, j+}(\sB_j, \sB_j) = \bSigma_{s, k}(\sB_j, \sB_j)$;
    \end{lsitemizeii}
\end{lsitemizei}

Now we have $\bSigma_{r, k}(\sB_i \cup \sB_j, \sB_i \cup \sB_j) = \bSigma_{s, k}(\sB_i \cup \sB_j, \sB_i \cup \sB_j)$. By
Property~\ref{prop:bb=bu}, we have $\bSigma_{r, k}(\sS_k\cup \sB_k, \sS_k\cup \sB_k) = \bSigma_{s, k}(\sS_k\cup \sB_k, \sS_k\cup
\sB_k)$. By induction, the theorem is proved. }

\begin{cor} \label{cor:final-conclusion}
For any $r$ and $s$ such that $\sC_i\in \sT_{r}^+$ and $\sC_i\in \sT_{s}^+$, we have:
\begin{equation*}
  \bSigma_{r, i+}(\sB_i, \sB_i) = \bSigma_{s, i+}(\sB_i, \sB_i)
\end{equation*}
\end{cor}
This corollary shows that the update processes for two target clusters $\sC_r$ and $\sC_s$ share their partial results so we can
reuse them in our algorithm. The proof is trivial given Theorem \ref{thm:individual-update} and Theorem
\ref{thm:synchronized-update}.

\bibliographystyle{elsarticle-num}
\bibliography{dissertation}

\begin{thebibliography}{10}
\expandafter\ifx\csname url\endcsname\relax
  \def\url#1{\texttt{#1}}\fi
\expandafter\ifx\csname urlprefix\endcsname\relax\def\urlprefix{URL }\fi
\expandafter\ifx\csname href\endcsname\relax
  \def\href#1#2{#2} \def\path#1{#1}\fi

\bibitem{datta97}
S.~Datta, Electronic Transport in Mesoscopic Systems, Cambridge University
  Press, 1997.

\bibitem{anant02}
A.~Svizhenko, M.~P. Anantram, T.~Govindan, B.~Biegel, R.~Venugopal,
  Two-dimensional quantum mechanical modeling of nanotransistors, Journal of
  Applied Physics 91~(4) (2002) 2343--2354.

\bibitem{ghosh2004molecular}
A.~Ghosh, P.~Damle, S.~Datta, A.~Nitzan, {Molecular Electronics}, Materials
  Research Society bulletin (2004) 391.

\bibitem{martinez2007self}
A.~Martinez, M.~Bescond, J.~Barker, A.~Svizhenko, M.~Anantram, C.~Millar,
  A.~Asenov, {A self-consistent full 3-D real-space NEGF simulator for studying
  nonperturbative effects in nano-MOSFETs}, IEEE Transactions on Electron
  Devices 54~(9) (2007) 2213--2222.

\bibitem{li08}
S.~Li, S.~Ahmed, G.~Klimeck, E.~Darve, Computing entries of the inverse of a
  sparse matrix using the {FIND} algorithm, Journal of Computational Physics
  227 (2008) 9408--9427.

\bibitem{darve04}
E.~Darve, S.~Li, Y.~Teslyar, Calculating transport properties of nanodevices,
  in: Proceedings of SPIE, Philadelphia, PA, Vol. 5593, 2004, pp. 452--463.

\bibitem{li07}
S.~Li, S.~Ahmed, E.~Darve, Fast inverse using nested dissection for {NEGF},
  Journal of Computational Electronics 6 (2007) 187--190.

\bibitem{lake1997single}
R.~Lake, G.~Klimeck, R.~Bowen, D.~Jovanovic, {Single and multiband modeling of
  quantum electron transport through layered semiconductor devices}, Journal of
  Applied Physics 81 (1997) 7845.

\bibitem{borm2003hierarchical}
S.~B\"orm, L.~Grasedyck, W.~Hackbusch, Hierarchical matrices, Lecture note,
  Citeseer 21 (2003) 1--124.

\bibitem{hackbusch99}
W.~Hackbusch, A sparse matrix arithmetic based on {$\cH$}-matrices. {P}art {I}:
  {I}ntroduction to {$\cH$}-matrices, Computing 62~(2) (1999) 89--108.

\bibitem{chandrasekaran06a}
S.~Chandrasekaran, P.~Dewilde, M.~Gu, T.~Pals, X.~Sun, A.~van~der Veen,
  D.~White, Some fast algorithms for sequentially semiseparable
  representations, Siam Journal on Matrix Analysis and Applications 27~(2)
  (2006) 341--364.

\bibitem{chandrasekaran06b}
S.~Chandrasekaran, M.~Gu, T.~Pals, A fast {ULV} decomposition solver for
  hierarchically semiseparable representations, SIAM Journal on Matrix Analysis
  and Applications 28 (2006) 603--622.

\bibitem{chandrasekaran08}
S.~Chandrasekaran, P.~Dewilde, M.~Gu, W.~Lyons, T.~Pals, A fast solver for
  {HSS} representations via sparse matrices, Siam Journal on Matrix Analysis
  and Applications 29~(1) (2008) 67--81.

\bibitem{mamaluy03}
D.~Mamaluy, M.~Sabathil, P.~Vogl, Efficient method for the calculation of
  ballistic quantum transport, Journal of Applied Physics 93~(8) (2003)
  4628--33.

\bibitem{khan2007quantum}
H.~Khan, D.~Mamaluy, D.~Vasileska, {Quantum transport simulation of
  experimentally fabricated nano-FinFET}, IEEE Transactions on Electron Devices
  54~(4) (2007) 784.

\bibitem{sabathil2004prediction}
M.~Sabathil, D.~Mamaluy, P.~Vogl, {Prediction of a realistic quantum logic gate
  using the contact block reduction method}, Semiconductor Science and
  Technology 19 (2004) S137--S138.

\bibitem{mamaluy2004contact}
D.~Mamaluy, A.~Mannargudi, D.~Vasileska, M.~Sabathil, P.~Vogl, {Contact block
  reduction method and its application to a 10 nm MOSFET device}, Semiconductor
  Science and Technology 19 (2004) S118--S121.

\bibitem{mamaluy2004electron}
D.~Mamaluy, A.~Mannargudi, D.~Vasileska, {Electron Density Calculation Using
  the Contact Block Reduction Method}, Journal of Computational Electronics
  3~(1) (2004) 45--50.

\bibitem{birner-ballistic}
S.~Birner, C.~Schindler, P.~Greck, M.~Sabathil, P.~Vogl, {Ballistic quantum
  transport using the contact block reduction (CBR) method}, Journal of
  Computational Electronics 8~(3) (2009) 267--286.

\bibitem{tang2010probing}
J.~Tang, Y.~Saad, {A probing method for computing the diagonal of the matrix
  inverse}, Tech. rep., Tech. Report umsi-2010-42, Minnesota Supercomputer
  Institute, University of Minnesota, Minneapolis, MN, 2009 (2010).

\bibitem{svizhenko02}
A.~Svizhenko, M.~P. Anantram, T.~R. Govindan, B.~Biegel, Two-dimensional
  quantum mechanical modeling of nanotransistors, Journal of Applied Physics
  91~(4) (2002) 2343--54.

\bibitem{erisman1975computing}
A.~Erisman, W.~Tinney, {On computing certain elements of the inverse of a
  sparse matrix}, Communications of the ACM 18~(3) (1975) 177--179.

\bibitem{petersen08}
D.~Petersen, H.~S{\o}rensen, P.~Hansen, S.~Skelboe, K.~Stokbro, Block
  tridiagonal matrix inversion and fast transmission calculations, J. Comput.
  Phys. 227 (2008) 3174--3190.

\bibitem{lin2009fast}
L.~Lin, J.~Lu, L.~Ying, R.~Car, {Fast algorithm for extracting the diagonal of
  the inverse matrix with application to the electronic structure analysis of
  metallic systems}, Communications in Mathematical Sciences 7~(3) (2009)
  755--777.

\bibitem{lin2011selinv}
L.~Lin, C.~Yang, J.~Meza, J.~Lu, L.~Ying, {SelInv---An Algorithm for Selected
  Inversion of a Sparse Symmetric Matrix}, ACM Transactions on Mathematical
  Software (TOMS) 37~(4) (2011) 40.

\bibitem{george73}
A.~George, Nested dissection of a regular finite-element mesh, SIAM Journal on
  Numerical Analysis 10~(2) (1973) 345--63.

\bibitem{liu1974solution}
J.~Liu, {The solution of mesh equations on a parallel computer}, Dep. of
  Computer Science, Univ., 1974.

\bibitem{george1983row}
A.~George, E.~Ng, {On row and column orderings for sparse least squares
  problems}, SIAM Journal on Numerical Analysis 20~(2) (1983) 326--344.

\bibitem{gilbert1987parallel}
J.~Gilbert, E.~Zmijewski, {A parallel graph partitioning algorithm for a
  message-passing multiprocessor}, International Journal of Parallel
  Programming 16~(6) (1987) 427--449.

\bibitem{pothen1990partitioning}
A.~Pothen, H.~Simon, K.~Liou, {Partitioning sparse matrices with eigenvectors
  of graphs}, SIAM Journal on Matrix Analysis and Applications 11~(3).

\bibitem{heath1991parallel}
M.~Heath, E.~Ng, B.~Peyton, {Parallel algorithms for sparse linear systems},
  SIAM review 33~(3) (1991) 420--460.

\bibitem{petersen09}
D.~Petersen, S.~Li, K.~Stokbro, H.~S{\o}rensen, P.~Hansen, S.~Skelboe,
  E.~Darve, {A hybrid method for the parallel computation of Green's
  functions}, Journal of Computational Physics 228~(14) (2009) 5020--5039.

\bibitem{golub96}
G.~H. Golub, C.~F. Van~Loan, Matrix Computations, Third Edition, Johns Hopkins
  University Press, 1996.

\bibitem{bunch71}
J.~R. Bunch, B.~N. Parlett, Direct methods for solving symmetric indefinite
  systems of linear equations, SIAM Journal on Numerical Analysis 8 (1971)
  639--55.

\bibitem{bunch77}
J.~R. Bunch, L.~Kaufman, Some stable methods for calculating inertia and
  solving symmetric linear systems, Mathematics of Computation 31 (1977)
  162--79.

\bibitem{trottenberg2009multigrid}
U.~Trottenberg, T.~Clees, {Multigrid software for industrial applications ---
  from MG00 to SAMG}, 100 Volumes of ‘Notes on Numerical Fluid Mechanics’
  (2009) 423--436.

\bibitem{brandt2000general}
A.~Brandt, {General highly accurate algebraic coarsening}, Electronic
  Transactions on Numerical Analysis 10~(1).

\bibitem{fattebert2000towards}
J.~Fattebert, J.~Bernholc, {Towards grid-based $O(N)$ density-functional theory
  methods: Optimized nonorthogonal orbitals and multigrid acceleration},
  Physical Review B 62~(3) (2000) 1713.

\bibitem{bernholc2008recent}
J.~Bernholc, M.~Hodak, W.~Lu, {Recent developments and applications of the
  real-space multigrid method}, Journal of Physics: Condensed Matter 20 (2008)
  294205.

\bibitem{feng2006nonlinear}
G.~Feng, T.~Beck, {Nonlinear multigrid eigenvalue solver utilizing
  nonorthogonal localized orbitals}, Physica Status Solidi (b) 243~(5) (2006)
  1054--1062.

\bibitem{wijesekera2007efficient}
N.~Wijesekera, G.~Feng, T.~Beck, {Efficient multiscale algorithms for solution
  of self-consistent eigenvalue problems in real space}, Physical Review B
  75~(11) (2007) 115101.

\end{thebibliography}

\end{document}